\input amstex
\documentstyle{amsppt}
\loadbold
\magnification = 1100
\mathsurround = 1 pt
\loadbold
\define\si{\smallskip\noindent}
\define\bi{\bigskip\noindent}
\define\pr{\text{ Pr\/}}
\define\ex{\text{\bf E\/}}

\define\e{\varepsilon}
\define\la{\lambda}
\define\a{\alpha}

\define\ga{\gamma}
\topmatter
\title
How frequently  is a system of $2$-linear Boolean equations solvable?
\endtitle
\rightheadtext{Random equations}
\author
Boris Pittel and Ji-A Yeum
\endauthor
\affil
Ohio State University
\endaffil
\address 
Ohio State University, Columbus, Ohio, USA
\endaddress
\email
bgp\@math.ohio-state.edu, yeum\@math.ohio-state.edu
\endemail
\thanks 
Pittel's research supported in part by NSF Grants  DMS-0406024, DMS-0805996
\endthanks
\abstract
We consider a random system of equations $x_i+x_j=b_{(i,j)}
(\text{mod }2)$, $(x_u\in \{0,1\},\, b_{(u,v)}=b_{(v,u)}\in\{0,1\})$,
with the pairs $(i,j)$ from $E$, a symmetric subset of $[n]\times [n]$. 
$E$ is chosen uniformly at random among all such subsets of a given cardinality $m$; 
alternatively $(i,j)\in E$ with a given probability $p$, independently of all
other pairs. Also, given $E$, $\pr\{b_{e}=0\}=\pr\{b_e=1\}$ for each $e\in E$, 
independently of all other $b_{e^\prime}$.
It is well known that, as $m$ passes through $n/2$ ($p$ passes through $1/n$, resp.),
the underlying random graph $G(n,\#\text{edges}=m)$, 
($G(n,\pr(\text{edge})=p)$, resp.) undergoes a rapid transition, from essentially
a forest of many small trees to a graph with one large, multicyclic, component in a sea of
small tree components. 
We should expect then that the solvability probability decreases precipitously in the
vicinity of $m\sim n/2$ ($p\sim 1/n$), and indeed this probability is of order
$(1-2m/n)^{1/4}$, for $m<n/2$ ($(1-pn)^{1/4}$, for $p<1/n$, resp.). We 
show that in a near-critical phase  $m=(n/2)(1+\la n^{-1/3})$ ($p=(1+\la n^{-1/3})/n$, resp.),
$\la=o(n^{1/12})$, the system is solvable with probability asymptotic to
$c(\la)n^{-1/12}$, for some explicit function $c(\la)>0$. Mike Molloy noticed that
the Boolean system with $b_e\equiv 1$ is solvable iff the underlying graph is
$2$-colorable, and asked whether this connection might be used to determine
an order of probability of $2$-colorability in the near-critical case. 
We answer Mike's question affirmatively and show that probability of
$2$-colorability  is $\lesssim 2^{-1/4}e^{1/8}c(\lambda)n^{-1/12}$, and 
asymptotic
to $2^{-1/4}e^{1/8}c(\la)n^{-1/12}$ at a critical phase $\la=O(1)$, and for $\la\to -\infty$.
(Submitted to Electronic Journal of Combinatorics on September 7, 2009.)
\endabstract
\endtopmatter
\document
{\bf 1. Introduction.\/} A system of $2$-linear equations over $GF(2)$ with $n$
Boolean variables $x_1,\dots,x_n\in \{0,1\}$ is
$$
x_i+x_j=b_{i,j}\,(\text{mod }2),\quad  b_{i,j}=b_{j,i}\in \{0,1\};\quad (i\neq j).\tag 1.1
$$
Here the unordered pairs $(i,j)$ correspond to the edge set of a given graph $G$
on the vertex set $[n]$. The system (1.1) certainly has a solution when $G$ is a tree. It 
can be obtained by picking an arbitrary $x_i \in \{0, 1\}$ at a ÒrootÓ $i$ and determining 
the other $x_j$ recursively along the paths leading away from the root. There is, 
of course, a twin solution $\bar x_j = 1 - x_j,\, j \in [n]$ . Suppose $G$ is not a tree, i.e. 
$\ell(G) := e(G) - v(G)\ge 0$. If $T$ is a tree spanning $G$, then each of additional edges 
$e_1 , \dots , e_{\ell(G)+1}$ forms, together with the edges of $T$, a single cycle $C_t$, 
$t \le \ell(G) + 1$. Obviously, a solution 
${x_j (T )}$ of a subsystem of (1.1) induced by the edges of $T$ is 
a solution of (1.1) provided that 
$$
b_{i,j }= x_i (T ) + x_j (T ), \,\,(i, j ) = e_1 , . . . , e_{\ell(G)+1} ; \tag 1.2
$$
equivalently 
$$
\sum_{e\in E(C_t)}b_e = 0\, (\text{mod } 2),\quad t = 1, . . . , \ell(G) + 1. \tag 1.3
$$
So, intuitively, the more edges $G$ has the less likely it is that the system (1.1) has 
a solution. We will denote the number of solutions by $S(G)$.

In this paper we consider solvability of a random system (1.1). Namely $G$ is 
either the Bernoulli random graph $G(n, p) = G(n,\pr(\text{edge}) = p)$, or the Erd\H os- 
R\'enyi random graph $G(n, m) = G(n, \# \text{ of edges} = m)$. Further, conditioned on the 
edge set $E(G(n, p))$ ($E(G(n, m)$ resp.), $b_e$'s are independent, and $\pr(b_e = 1) = 
\hat p $,  for all $e$. We focus on  $\hat p = 1/2$ and 
$\hat p = 1$. $\hat p=1/2$ is the case when $b_e$'s are ``absolutely random''.
For $\hat p=1$, $b_e$'s are all ones.  Mike Molloy [17], who brought this case to
our attention, noticed that here (1.1) has a solution iff the underlying graph is bipartite, 
2-colorable in other words.

 It is well known that, as m passes through $n/2$ ($p$ passes through $1/n$, resp.), the 
underlying random graph $G(n, m)$, ($G(n, p)$, resp.) undergoes a rapid transition, 
from essentially a forest of many small trees to a graph with one large, multicyclic, 
component in a sea of small tree components. Bollob\'as [4], [5] discovered that, for 
$G(n, m)$, the phase transition window is within $[m_1,m_2 ]$ , where 
$$
m_{1,2} = n /2 \pm \la n^{2/3 }, \quad \la = \Theta(\ln^{1/2} n). 
$$
\L uczak [14] was able to show that the window is precisely $[m_1,m_2 ]$ with $\la
\to\infty$ however slowly. (See  \L uczak et al [16], Pittel [19] for the distributional
results on the critical graphs $G(n,m)$ and $G(n,p)$.) We should expect then that the solvability 
probability decreases 
precipitously for $m$ close to $n/2$ ($p$ close to $1/n$ resp.). Indeed, for a {\it multi\/}graph
version of $G(n,m)$, Kolchin [13] proved that this probability 
is asymptotic to 
$$
\frac{(1- \ga) ^{1/4}}{(1-(1-2\hat p)\ga)^{1/4}},\quad \ga:=\frac{2m}{n},\tag 1.4
$$
if  $\limsup\ga<1$. 
See Creignon and Daud\'e [9] for a similar result. Using the results from Pittel 
[19], we show (see Appendix)
that for the random graphs $G(n,\ga n/2)$ and $G(n,p=\ga/n)$, with $\limsup\ga<1$,
the corresponding probability is asymptotic to
$$
\frac{(1- \ga) ^{1/4}}{(1-(1-2\hat p)\ga)^{1/4}}\exp\left[\frac{\ga}{2}\hat p+
\frac{\ga^2}{2}\hat p(1-\hat p)\right].\tag 1.5
$$
The relations (1.4), (1.5) make it plausible that, in the nearcritical
phase $|m-n/2|=O(n^{2/3})$, the solvability probability is of order $n^{-1/12}$. Our 
goal
is to confirm, rigorously, this conjecture. 
\si

To formulate our main result, we need some notations. Let $\{f_r\}_{r\ge 0}$ be a sequence
defined by an implicit recurrence
$$
f_0=1,\quad \sum_{k=0}^r f_kf_{r-k}=\varepsilon_r,\quad \varepsilon_r:=\frac{(6r)!}{2^{5r}3^{2r}(3r)!(2r)!}.
\tag 1.6
$$
Equivalently, the formal series $\sum_rx^rf_r$, $\sum_rx^r\varepsilon_r$ (divergent for all 
$x\neq 0$)
satisfy
$$
\left(\sum_rx^rf_r\right)^2=\sum_rx^r\varepsilon_r.\tag 1.7
$$
It is not difficult to show that
$$
\frac{\varepsilon_r}{2}\left(1-\frac{1}{r}\right)\le f_r\le \frac{\varepsilon_r}{2},\quad r>0.
\tag 1.8
$$
For $y,\la\in \Bbb R$, let $A(y,\la)$ denote the sum of a convergent series,
$$
A(y,\la)=\frac{e^{-\la^3/6}}{3^{(y+1)/3}}\sum_{k\ge 0}
\frac{\left(\frac{1}{2}3^{2/3}\la\right)^k}
{k!\Gamma [(y+1-2k)/3]}.\tag 1.9
$$
We will write $B_n\sim C_n$ if $\lim_{n\to\infty}B_n/C_n= 1$, and $B_n\lesssim C_n$ if
$\limsup_n B_n/C_n\le 1$.
Let $S_n$ denote the random number of solutions
of (1.1) with the underlying graph being either $G(n,m)$ or $G(n,p)$, i. e.
$S_n=S(G(n,m))$ or $S_n=S(G(n,p))$, and 
the (conditional) probability of $b_e=1$ for $e\in E(G(n,m))$ ($e\in E(G(n,p))$ resp.)
being equal $\hat p$.
\proclaim{Theorem 1.1} {\bf (i)\/} Let $\hat p=1/2$. Suppose that
$$
m=\frac{n}{2}\,(1+\la n^{-1/3}),\quad p=\frac{1+\la n^{-1/3}}{n},\quad\quad |\la| =o(n^{1/12}).\tag 1.10
$$
Then, for both $G(n,m)$ and $G(n,p)$,
$$
\pr(S_n>0)\sim\,n^{-1/12}c(\la),\tag 1.11
$$
where
$$
c(\la):=\left\{\alignedat2
&e^{3/8}(2\pi)^{1/2}\sum_{r\ge 0}
\frac{f_r}{2^r}\,A(0.25+3r,\la),\quad&&\la\in (-\infty,\infty);\\
&e^{3/8}|\la|^{1/4},\quad&&\la\to-\infty;\\
&\frac{e^{3/8}}{4\cdot 3^{3/4}}\,\la^{1/4}\exp(-10\la^3/81),\quad&&\la\to\infty.
\endalignedat\right.\tag 1.12
$$
{\bf (ii)\/} Let $\hat p=1$. Then, with $c(\la)$ replaced by
$c_1(\la):=2^{-1/4}e^{1/8}c(\la)$, (1.9) holds for both $G(n,m)$ and $G(n,p)$ 
if either $\lambda=O(1)$, or
$\la\to-\infty$, $|\la|=o(n^{1/12})$. For $\la\to\infty$, $\la=o(n^{1/12})$,
$$
\pr(S_n>0)\lesssim n^{-1/12}c_1(\la).
$$
\endproclaim
\si
 
{\bf Notes.\/} 1. For $G(n,m)$ with $\la\to-\infty$, and $\hat p=1/2$, our result 
blends, qualitatively, with 
the estimate (1.4)  from [13] and [9]  for a subcritical {\it multi\/}graph,  and  becomes 
the
estimate (1.5) for the subcritical graphs $G(n,m)$ and $G(n,p)$. 
\si
2. The part {\bf (ii)\/} answers Molloy's question: the critical graph $G(n,m)$
($G(n,p)$ resp.) is bichromatic (bipartite) with probability
$\sim c_1(\la)n^{-1/12}$. 
\si

Very interestingly, the largest bipartite subgraph of the critical $G(n,p)$
can be found in expected time $O(n)$, see Coppersmith et al [8], Scott and Sorkin [21] and
references therein. The case $\la\to\infty$ of {\bf (ii)\/} strongly
suggests that the supercritical graph $G(n,p=c/n)$, ($G(n,m=cn/2)$ resp.),
i. e. with $\liminf c>1$, is bichromatic with exponentially small
probability. In [8] this exponential smallness  was established for
the conditional probability, given that the random graph has a giant component. 
\si

Here is a technical reason why, for $\la=O(1)$ at least,  the asymptotic probability of 
$2$-colorability is 
the asymptotic solvability probability for (1.1) with $\hat p=1/2$
{\it times\/} $2^{-1/4}e^{1/8}$. Let $C_{\ell}(x)$ ($C_{\ell}^e(x)$ 
resp.)
denote the exponential generating functions of connected graphs $G$ (graphs $G$ without
odd cycles resp.) with excess  $e(G)-v(G)=\ell\ge 0$. It turns out that, for $|x|<e^{-1}$
 (convergence radius of $C_{\ell}(x)$, $C_{\ell}^e(x)$), 
and $x\to e^{-1}$,
$$
C_{\ell}^e(x)\left\{\alignedat2
&\sim \frac{1}{2^{\ell+1}}C_{\ell}(x),\quad&&\ell>0,\\
&=\frac{1}{2}C_0(x)+\ln\bigl(2^{-1/4}e^{1/8}\bigr)+o(1),\quad&&\ell=0.
\endalignedat\right.
$$  
Asymptotically, within the factor $e^{\ln\bigl(2^{-1/4}e^{1/8}\bigr)}$, this 
reduces 
the problem to that for $\hat p=1/2$. Based on (1.5),
we conjecture that generally, for
$\hat p\in (0,1]$, and the critical $p$, $\pr(S_n>0)$ is
that probability for $\hat p=1/2$ times
$$
(2\hat p)^{-1/4}\exp\left[-\frac{(1-\hat p)^2}{2}+\frac{1}{8}\right].
\tag 1.11
$$ 
(For $\hat p=0$, $\pr(S_n>0)=1$ 
obviously.)
\si

3. While working on this project, we became aware of a recent paper
[10] by Daud\'e and Ravelomanana. They  studied  a close but different case, when
a system of $m$ equations is chosen uniformly at random among {\it all\/} $n(n-1)$ equations 
of the form (1.1).
In particular, it is possible to have  pairs of clearly contradictory equations, $x_i+x_j=0$ and 
$x_i+x_j=1$. For $m=O(n)$ the probability that none of these simplest contradictions occurs
is bounded away from zero. So, intuitively, the system they studied is close to ours with
$G=G(n,m)$ and $\hat p=1/2$. Our asymptotic formula (1.9), with two first equations in
(1.10), in this case is similar to 
Daud\'e-Ravelomanana's  main theorem, 
but there are some puzzling differences. The exponent series in their equation (2) is certainly
misplaced; their claim  does not contain our sequence $\{f_r\}$. 
\si
As far as we can judge by a proof outline in [10], our argument is quite different. Still
like [10],
our analysis is based on the generating functions of sparse graphs
discovered, to a great extent, by Wright [23],  [24].  We gratefully credit Daud\'e and
Ravelomanana for stressing importance of Wright's  bounds for the
generating function $C_{\ell}(x)$. These bounds  play a substantial role in our 
argument as well. 
\si

4. We should mention a large body of work on a related, better known, $2-SAT$ problem, see for 
instance
Bollob\'as et al [6], and references therein. It is a problem of existence of a 
truth-satisfying assignment for the variables in the conjunction of $m$ random disjunctive clauses 
of 
a form $x_i\lor x_j$, ($i,j\in [n]$). It is well known, Chv\'atal and Reed [7],  that the 
existence threshold 
is $m/n=1$. It was proved in [6]  that the phase transition window is $[m_1,m_2]$, with
$$
m_{1,2}=n\pm \la\, n^{2/3},\quad |\la|\to\infty\text{  however slowly},
$$
and that the solvability probability is bounded away from both $0$ and $1$ iff
$m=n+O(n^{2/3})$.
\si

5. A natural extension of the system (1.1) is a system of $k$-linear equations
$$
\sum_{i\in e}x_i =b_e\,(\text{mod }2), \tag 1.12
$$
where $e$ runs over a set $E$ of (hyper)edges of a $k$-uniform hypergraph $G$, $k\ge 2$, on the vertex
set $[n]$, Kolchin [13]. Suppose $G$ is chosen uniformly at random among all $k$-uniform graphs with
a given number $m$ of edges, and, given $G$, the $b_e$s are independent Bernoullis. It will be
interesting to study, for $k>2$, the limiting solvability probability as a function of $m/n$. See
[13] for some thought-provoking results on the behavior of the number of
hypercycles in this random hypergraph. 
\bi

The paper is organized as follows. In the section 2 we work on the $G(n,p)$ and
$\hat p=1/2$ case. 
\si
Specifically in the (sub)section 2.1 we express the
solvability probability, $\pr(S_n>0)$, and its truncated version, as a coefficient by $x^n$
in a power series based on the generating functions of the sparsely edges
(connected) graphs. We also establish positive correlation between
solvability and boundedness of a maximal ``excess'', and determine a proper
truncation of the latter dependent upon the behavior of $\la$. In the section
(2.2) we provide a necessary information about the generating functions
and their truncated versions involved in the formula and the bounds for $\pr(S_n>0)$.
In the section 2.3 we apply complex analysis techniques to the ``coefficient
by $x^n\,$'' formulas and obtain a sharp
asymptotic estimate for $\pr(S_n>0)$ for $|\la|=o(n^{1/12})$.

In the section 3 we transfer the results of the section 2 to the $G(n,m)$ and
$\hat p=1/2$ case .

In the section 4 we establish the counterparts of the results from the sections
2,3 for $G(n,p)$, $G(n,m)$ with $\hat p=1$. An enumerative ingredient of
the argument is an analogue of Wright's formulas for the generating
functions of the connected graphs without odd cycles.

In Appendix we prove some auxilliary technical results, and an asymptotic
formula for $\pr(S_n>0)$ in the subcritical case, i. e. when
the average vertex degree is less than, and bounded away from $1$.

\bi
{\bf 2. Solvability probability: $G(n,p)$ and $\hat p=1/2$.\/}
\bi 
{\bf 2.1. Representing bounds for $\pr(S_n>0)$ as a coefficient of  $x^n$ in a power series.\/} 
\bi
Our first step is to compute the probability of the event $\{S_n>0\}$, conditioned on $G(n,p)$.
Given a graph $G=(V(G),E(G))$, we denote $v(G)=|V(G)|$,
$e(G)=|E(G)|$.
  
\proclaim{Lemma 2.1.1} Given a graph $G$ on $[n]$, 
let $c(G)$ denote the total number of its components $H_i$. Then
$$
\align
\pr(S_n>0\,|\,G(n,p)=G)=&\prod_{i=1}^{c(G)}\left(\frac{1}{2}\right)^{e(H_i)-(v(H_i)-1)}\\
=&\left(\frac{1}{2}\right)^{X(G)},\quad X(G):=e(G)-n+c(G).
\endalign
$$
Consequently
$$
\pr(S_n>0)=\ex\left[\left(\frac{1}{2}\right)^{X(G(n,p))}\right].
$$
\endproclaim

{\bf Proof of Lemma 2.1.1.\/} Recall that, conditioned on $G(n,p)$, the edge variables
$b_e$ are mutually independent. So it is suffices to show that a system (1.1) for a connected graph
$H$, with independent $b_e$, $e\in E(H)$, such that $\pr(b_e=1)=1/2$, is solvable with
probability $(1/2)^{\ell +1}$, where $\ell=e(H)-v(H)$.
 
Let $T$ be a tree spanning $H$. Let $\bold x(T):=\{x_i(T)\}_{i\in V(H)}$ be the solution
of the subsystem of (1.1) corresponding to $v(H)-1$
edges of $T$, with $x_{i_0}=1$ say, for a specified ``root'' $i_0$. $\bold x(T)$
is a solution of the whole system (1.1) iff 
$$
b_e=x_i(T)+x_j(T),\quad ((i,j)=e), \tag 2.1.1
$$
for each of $e(H)-(v(H)-1)=\ell+1$ edges $e\in E(H)\setminus E(T)$. By independence of $b_e$'s,
the probability that, conditioned on  $\{b_e\}_{e\in E(T)}$,  the constraints
(2.1.1) are met is $(1/2)^{\ell+1}$,. (It is crucial that $\pr(b_e=0)=\pr(b_e=1)=1/2$.)
Hence the unconditional solvability probability for the system (1.1) with the underlying
graph $H$ is $(1/2)^{\ell+1}$ as well. \qed
\si

{\bf Note.\/} For a cycle $C\subseteq H$, let $b_C=\sum_{e\in E(C)}b_e$. The conditions
(2.1.1) are equivalent to $b_C$ being even for the $\ell+1$ cycles, each formed by adding to $T$ an
edge in $E(H)\setminus E(T)$. Adding the equations (1.1) over the edges of
{\it any\/} cycle $C\subseteq H$, we see that necessarily $b_C$ is even too.  
Thus our proof effectively shows that 
$$
\pr\left\{\bigcap_{C\subseteq H}\{b_C\text{ is even}\}\right\}=\left(\frac{1}{2}\right)
^{\ell(H)+1}.
$$
\bi

Using Lemma 2.1.1, we express $P(S(n,p)>0)$ as the coefficient by $x^n$ in a {\it formal\/}
power series.
To formulate the result, introduce $C_{\ell}(x)$, the exponential generating function
of a sequence $\{C(k,k+\ell)\}_{k\ge 1}$, where $C(k,k+\ell)$ is the total number of
connected graphs $H$ on $[k]$ with {\it excess\/} $e(H)-v(H)=\ell$. Of course, $C(k,k+\ell)=0$
unless $-1\le \ell\le \binom{k}{2}-k$.
\proclaim{Lemma 2.1.2} 
$$
\align
\pr(S_n>0)=&\, N(n,p)\,[x^n]\,\exp\left[\frac{1}{2}\sum_{\ell\ge -1}
\left(\frac{p}{2q}\right)^{\ell}C_{\ell}(x)\right],\tag 2.1.2\\
N(n,p):=&\,n!\,q^{n^2/2}\left(\frac{p}{q^{3/2}}\right)^n.\tag 2.1.3
\endalign
$$
\endproclaim

{\bf Proof of Lemma 2.1.2.\/} The proof mimicks derivation of the ``coefficient-of $x^n$-
expression'' for the largest component size distribution in [19].
\si
Given $\boldsymbol\alpha=\{\alpha_{k,\ell}\}$, such 
that $\sum_{k,\ell}k\alpha_{k,\ell}=n$,
let $P_n(\boldsymbol\alpha)$ denote the probability that $G(n,p)$ has $\alpha_{k,\ell}$ components
$H$ with $v(H)=k$ and $e(H)-v(H)=\ell$. To compute 
$P_n(\boldsymbol\alpha)$, we observe that there are
$$
\frac{n!}{\prod\limits_{k,\ell}(k!)^{\alpha_{k,\ell}}\alpha_{k,\ell}!}
$$
ways to partition $[n]$ into $\sum_{k,\ell}\alpha_{k,\ell}$ subsets, with $\alpha_{k,\ell}$
subsets of cardinality $k$ and ``type'' $\ell$. For each such partition, there are
$$
\prod_{k,\ell}[C(k,k+\ell)]^{\alpha_{k,\ell}}
$$
ways to build $\alpha_{k,\ell}$ connected graphs $H$ on the corresponding $\alpha_{k,\ell}$
subsets, with $v(H)=k$, $e(H)-v(H)=\ell$.
The probability that these graphs are induced subgraphs of $G(n,p)$ is
$$
\prod_{k,\ell}\left[p^{k+\ell}q^{\binom{k}{2}-(k+\ell)}\right]^{\alpha_{k,\ell}}=
\left(\frac{p}{q^{3/2}}\right)^n\prod_{k,\ell}\left[\left(\frac{p}{q}\right)^{\ell}q^{k^2/2}
\right]^{\alpha_{k,\ell}},
$$
as $\sum_{k,\ell}k\,\alpha_{k,\ell}=n$. The probability that no two vertices from two different subsets are joined by an edge
in $G(n,p)$ is $q^r$, where $r$ is the total number of all such pairs, i. e.
$$
\align
r=&\sum_{k,\ell}k^2\binom{\alpha_{k,\ell}}{2}+
\frac{1}{2}\sum_{(k_1,\ell_1)\neq (k_2,\ell_2)}k_1k_2\alpha_{k_1,\ell_1}\alpha_{k_2,\ell_2}\\
=&-\frac{1}{2}\sum_{k,\ell}k^2\alpha_{k,\ell}+\frac{1}{2}\left(\sum_{k,\ell}k\,\alpha_{k,\ell}
\right)^2\\
=&-\frac{1}{2}\sum_{k,\ell}k^2\alpha_{k,\ell}+\frac{n^2}{2}.
\endalign
$$
Multiplying the pieces,
$$
P_n(\boldsymbol\alpha)=N(n,p)\prod_{k,\ell}\frac{1}{\alpha_{k,\ell}!}
\left[\frac{(p/q)^{\ell}C(k,k+\ell)}{k!}\right]^{\alpha_{k,\ell}}.
$$
So, using Lemma 2.1.1,
$$
\pr(S_n>0)=N(n,p)\sum_{\boldsymbol\alpha}\prod_{k,\ell}\frac{1}{\alpha_{k,\ell}!}
\left[\frac{(1/2)^{\ell+1}(p/q)^{\ell}C(k,k+\ell)}{k!}\right]^{\alpha_{k,\ell}}.\tag 2.1.4
$$

Notice that dropping factors $(1/2)^{\ell+1}$ on the right, we get $1$ instead of
$\pr(S_n>0)$ on the left, i.e.
$$
1=N(n,p)\sum_{\boldsymbol\alpha}\prod_{k,\ell}\frac{1}{\alpha_{k,\ell}!}
\left[\frac{(p/q)^{\ell}C(k,k+\ell)}{k!}\right]^{\alpha_{k,\ell}}.\tag 2.1.5
$$

So, multiplying both sides of (2.1.4) by $\frac{x^n}{N(n,p)}$ and summing over $n\ge 0$,
$$
\aligned
\sum_n x^n\,\frac{\pr(S_n>0)}{N(n,p)}=&\sum_{\sum\limits_{k,\ell}k\alpha_{k,\ell}<\infty}
\,\,\prod_{k,\ell}\frac{x^{k\alpha_{k,\ell}}}{\alpha_{k,\ell}!}
\left[\frac{(1/2)^{\ell+1}(p/q)^{\ell}C(k,k+\ell)}{k!}\right]^{\alpha_{k,\ell}}\\
=&\exp\left[\frac{1}{2}\sum_{\ell}(p/2q)^{\ell}\sum_k\frac{C(k,k+\ell)x^k}{k!}\right]\\
=&\exp\left[\frac{1}{2}\sum_{\ell}(p/2q)^{\ell}C_{\ell}(x)\right].
\endaligned\tag 2.1.6
$$
We hasten to add that the series on the right, whence the one on the left, converges for
$x=0$ only. Indeed, using (2.1.5) instead of (2.1.4),
$$
\exp\left[\sum_{\ell}(p/q)^{\ell}C_{\ell}(x)\right]=\sum_n x^n\,\frac{1}{N(n,p)}=
\sum_n \frac{\left(\frac{xq^{3/2}}{p}\right)^n}{n!\,q^{n^2/2}}=\infty,\tag 2.1.7
$$
for each $x>0$. Therefore, setting $p/2q=p_1/q_1$, ($q_1=1-p_1$),
$$
\sum_{\ell}(p/2q)^{\ell}C_{\ell}(x)=\sum_{\ell}(p_1/q_1)^{\ell}C_{\ell}(x)=\infty,\quad
\forall x>0,
$$
as well.\qed
\si

{\bf Note.\/} Setting $p/q=w$, $x=yw$, in (2.1.7), so that $p=w/(w+1)$, $q=1/(w+1)$, we obtain
a well known (exponential) identity, e. g. Janson et al [12],  
$$
\exp\left[\sum_{\ell\ge -1}w^{\ell}C_{\ell}(yw)\right]=\sum_{n\ge 0}\frac{y^n}{n!}
(w+1)^{\binom{n}{2}};
$$
the right expression (the left exponent resp.) is a bivariate generating function
for graphs (connected graphs resp.) $G$ enumerated by $v(G)$ and $e(G)$.
Here is a similar identity involving generating functions of connected
graphs $G$ with a fixed positive excess, 
$$ 
\exp\left[\sum_{\ell\ge 1}w^{\ell}C_{\ell}(x)\right]=\sum_{r\ge 0}w^{r}E_{r}(x),
\tag 2.1.8
$$
where $E_0(x)\equiv 1$, and, for $\ell\ge 1$, $E_{\ell}(x)$ is the exponential generating
function of graphs $G$ without tree components and unicyclic components, that have excess 
$\ell(G)=e(G)-v(G)=\ell$, see [12]. In the light of Lemma 2.1.2, we will
need an expansion
$$ 
\exp\left[\frac{1}{2}\sum_{\ell\ge 1}w^{\ell}C_{\ell}(x)\right]=\sum_{r\ge 0}w^{r}F_{r}(x).
\tag 2.1.9
$$
Like $E_r(x)$, each power series $F_r(x)$ has nonnegative coefficients, and converges
for $|x|<e^{-1}$.  
\bi

By Lemma 2.1.2 and (2.1.8),
$$
\aligned
\pr(S_n>0)=&N(n,p)\sum_{r\ge 0}\left(\frac{p}{2q}\right)^r[x^n]\left\{
e^{H(x)}F_r(x)\right\};\\
H(x):=&\frac{q}{p}C_{-1}(x)+\frac{1}{2}C_0(x).
\endaligned\tag 2.1.10
$$
Interchange of $[x^n]$ and the summation is justifiable as each of the functions
on the right has a power series expansion with only nonnegative coefficients. 
That is, divergence of $\sum_{\ell}(p/2q)^{\ell}C_{\ell}(x)$ in (2.1.6) does not impede
evaluation of $\pr(S_n>0)$. Indirectly though this divergence does make it difficult, if
possible at all, to obtain a sufficiently sharp estimate of the terms in the above
sum for $r$ going to $\infty$ with $n$, needed to derive an asymptotic formula
for that probability. Thus
we need to truncate, one way or another, the divergent series on the right in (2.1.6). 
One of the properties of $C_{\ell}(x)$ discovered
by Wright [23] is that each of these series converges (diverges) for $|x|<e^{-1}$
(for $|x|>e^{-1}$ resp.). So, picking $L\ge 0$, and restricting summation range 
to $\ell\in [-1,L]$, we definitely get a series convergent for $|x|<e^{-1}$. What is then
a counterpart of $\pr(S_n>0)$? Perusing the proof of Lemma 2.1.2, we easily see the
answer.
\si

Let $G$ be a graph with components
$H_1,H_2,\dots$. Define $\Cal E(G)$, a maximum excess of $G$, by
$$
\Cal E(G)=\max_i[e(H_i)-v(H_i)].
$$
It can be easily seen that $\Cal E(G)$ is monotone increasing, i. e. $\Cal E(G^\prime)
\le\Cal E(G^{\prime\prime})$ if $G^\prime\subseteq G^{\prime\prime}$. Let
$\Cal E_n=\Cal E(G(n,p))$.
\proclaim{Lemma 2.1.3}
$$
\pr(S_n>0,\,\Cal E_n\le L)=\, N(n,p)\,[x^n]\,\exp\left[\frac{1}{2}\sum_{\ell=-1}^L
\left(\frac{p}{2q}\right)^{\ell}C_{\ell}(x)\right],\tag 2.1.11
$$
\endproclaim
\noindent
The proof of (2.1.11) is an obvious modification of that for (2.1.2).
\bi
If, using (2.1.11), we are able to estimate $\pr(S_n>0,\,\Cal E_n\le L)$, then evidently
we will get a {\it lower\/} bound of $\pr(S_n>0)$, via
$$
\pr(S_n>0)\ge \pr(S_n>0,\,\Cal E_n\le L).\tag 2.1.12
$$
Crucially, the events $\{S_n>0\}$ and $\{\Cal E_n\le L\}$ are {\it positively\/}
correlated.
\proclaim{Lemma 2.1.4}
$$
\pr(S_n>0)\le\frac{\pr(S_n>0,\,\Cal E_n\le L)}{\pr(\Cal E_n\le L)}.\tag 2.1.13
$$
\endproclaim
\si

{\bf Note.\/} The upshot of (2.1.12)-(2.1.13) is that 
$$
\pr(S_n>0)\sim\pr(S_n>0,\Cal E_n\le L),
$$
provided that $L=L(n)$ is just large enough to guarantee 
that $\pr(\Cal E_n\le L)\to 1$. 
\si

{\bf Proof of Lemma 2.1.4.\/} By Lemma 2.1.1,
$$
\pr(S_n>0,\,\Cal E_n\le L)=\ex\left[\left(\frac{1}{2}\right)^{X(G(n,p))}
\bold 1_{\{\Cal E(G(n,p))\le L\}}\right],
$$
where $X(G)=e(G)-n+c(G)$. Notice that $(1/2)^{X(G)}$ is monotone decreasing. Indeed, if a graph
$G_2$ is obtained by adding one edge to a graph $G_1$, then 
$$
e(G_2)=e(G_1)+1,\quad c(G_2)\in \{c(G_1)-1,c(G_1)\},
$$
so that $X(G_2)\ge X(G_1)$. Hence, using induction on $e(G_2)-e(G_1)$,
$$
G_1\subseteq G_2\Longrightarrow X(G_2)\ge X(G_1).
$$
Furthermore $\bold 1_{\{\Cal E(G)\le L\}}$ is also monotone decreasing. (For
$e\notin E(G)$, if $e$ joins two vertices from the same component of $G$
then $\Cal E(G+e)\ge \Cal E(G)$ obviously. If $e$ joins two components, $H_1$
and $H_2$ of $G$, then the resulting component has an excess more than or equal to
$\max\{\Cal E(H_1),\Cal E(H_2)\}$, with equality when one of two components is a tree.)
\si

Now notice
that each $G$ on $[n]$ is essentially a $\binom{n}{2}$-long tuple
$\boldsymbol\delta$ of $\{0,1\}$-valued
variables $\delta_{(i,j)}$, $\delta_{(i,j)}=1$ meaning that $(i,j)\in E(G)$.
So, a graph function $f(G)$ can be unambigiously written as $f(\boldsymbol\delta)$.
Importantly, a monotone decreasing (increasing) graph function is a monotone
decreasing (increasing) function of the code $\boldsymbol\delta$. For the random
graph $G(n,p)$, the components of $\boldsymbol\delta$ are independent random variavbles. 
According to an FKG-type inequality, see Grimmett and Stirzaker [11] for instance,
for any two decreasing (two increasing) functions $f(\bold Y)$, $g(\bold Y)$ of a vector
$\bold Y$ with independent components,
$$
\ex[f(\bold Y)g(\bold Y)]\ge \ex[f(\bold Y)]\,\ex[g(\bold Y)].
$$
Applying this inequality to $(1/2)^{X(\boldsymbol\delta)}
\bold 1_{\{\Cal E(\boldsymbol\delta)\le L\}}$, we obtain
$$
\align
\pr(S_n>0,\,\Cal E_n\le L)\ge& \ex\left[\left(\frac{1}{2}\right)^{X(G(n,p))}\right]
\ex\left[\bold 1_{\{\Cal E(G(n,p))\le L\}}\right]\\
=&\pr(S_n>0)\pr(\Cal E_n\le L).
\endalign
$$
\qed
\bi

Thus our next step is to determine how large $\Cal E(G(n,p))$ is typically,
if
$$
p=\frac{1+\la n^{-1/3}}{n},\quad \la=o(n^{1/3}). \tag 2.1.14
$$
For $p=c/n$, $c<1$, it was shown in Pittel [19] that 
$$
\lim\pr(G(n,p)\text{ does not have a cycle})=(1-c)^{1/2}\exp(c/2+c^2/4).
$$
From this result and monotonicity of $\Cal E(G)$, it follows that, for $p$ in 
(2.1.14), 
$$
\lim\pr(\Cal E(G(n,p))\ge 0)=1.
$$
If $\la\to -\infty$, then we also have
$$
\lim\pr(\Cal E(G(n,p))>0)=0,\tag 2.1.15
$$
that is $\Cal E(G(n,p))\le 0$ with high probability (whp). (The proof of (2.1.15)
mimicks \L uczak's 
proof [14] of an analogous property of $G(n,m)$, with $n^{-2/3}(n/2-m)\to\infty$.)
\si
Furthermore, by Theorem 1 in [16],  and monotonicity of $\Cal E(G(n,p))$, it follows that $\Cal E(G(n,p))$
is bounded in probability  (is $O_{P}(1)$, in short), if $\limsup\la<\infty$.
\si
Finally, suppose that $\la\to\infty$. Let $\Cal L(G(n,m))$ denote the total excess of the number of
edges over the number of vertices in the {\it complex\/} components of $G(n,m)$, i. e.
the components that are neither trees nor unicyclic. According to a limit theorem for 
$\Cal L(G(n,m=(n/2)(1+\la n^{-1/3})))$ from [12], $\Cal L(G(n,m))/\la^3\to 2/3$, in probability. 
According to \L uczak [14], whp $G(n,m)$ has exactly one complex component. So whp
$\Cal E(G(n,m))=\Cal L(G(n,m))$, i. e. $\Cal E(G(n,m))/\la^3\to 2/3$ in probability, as well.

Now, if 
$$
m^\prime=Np+O\left(\sqrt{Npq}\right),\quad N:=\binom{n}{2},
$$
then
$$
m^\prime=\frac{n}{2}(1+\la^\prime n^{-1/3}),\quad \la^\prime:=\la\bigl(1+O(n^{-1/6})\bigr).
$$
Therefore, in probability,
$$
\frac{\Cal E(G(n,m^\prime))}{\la^3}\to \frac{2}{3},
$$
as well.  From a general ``transfer principle'' ([5], [15])  it follows then that
$$
\frac{\Cal E(G(n,p))}{\la^3}\to \frac{2}{3},
$$
in probability, too.
\bi

This  discussion justifies the following choice of $L$:
$$
L=\left\{\alignedat2
&0,\quad&&\text{if }\lim\la=-\infty,\\
&u\to\infty\text{ however slowly},\quad&&\text{if }\la=O(1),\\
&\la^3,\quad&&\text{if }\la\to\infty,\,\la=o(n^{1/12}).\endalignedat\right.
\tag 2.1.16
$$
\bi
{\bf 2.2 Generating functions.\/}
\bi

First, some basic facts about the generating functions $C_{\ell}(x)$ and $E_{\ell}(x)$.  Introduce
a tree function $T(x)$, the exponential generating function of $\{k^{k-1}\}$, the counts
of rooted trees on $[k]$, $k\ge 1$. It is well known that the
series 
$$
T(x)=\sum_{k\ge 1}\frac{x^k }{k!}k^{k-1}
$$
has convergence radius $e^{-1}$, and that  
$$
T(x)=xe^{T(x)},\quad |x|\le e^{-1};
$$
in particular, $T(e^{-1})=1$. (This last fact has a probabilistic explanation:
$\{\frac{k^{k-1}}{e^kk!}\}$ is the distribution of a total progeny in a branching
process with an immediate family size being Poisson ($1$) distributed.) $T(x)$ is a building
block for all $C_{\ell}(x)$. Namely, (Moon [18], Wright [23], Bagaev [1] resp.),
$$
\align
C_{-1}(x)=&\, T(x)-\frac{1}{2}T^2(x),\tag 2.2.1\\
C_0(x)=&\,\frac{1}{2}\left[\ln\frac{1}{1-T(x)}-T(x)-\frac{1}{2}T^2(x)\right],\tag 2.2.2\\
C_1(x)=&\,\frac{T^4(x)(6-T(x))}{24(1-T(x))^3},
\endalign
$$
and ultimately, for all $\ell>0$,
$$
C_{\ell}(x)=\sum_{d=0}^{3\ell+2}\frac{c_{\ell,d}}{(1-T(x))^{3\ell-d}},\tag 2.2.3
$$
Wright [23]. Needless to say, $|x|<e^{-1}$ in all the formulas. One should rightfully
anticipate though that the behaviour of $C_{\ell}(x)$ for $x$'s close to $e^{-1}$ is going to
determine an asymptotic behaviour of $\pr(S_n>0,\Cal E_n\le L)$. And so the $(d=0)$-term
in (2.2.3) might well be the only term we would need eventually. In this context, it is
remarkable that in a follow-up
paper [24] Wright was able to show that 
$$
\multline
\frac{c_{\ell}}{(1-T(x))^{3\ell}}-\frac{d_{\ell}}{(1-T(x))^{3\ell-1}}\le_c
C_{\ell}(x)\\
\le_c \frac{c_{\ell}}{(1-T(x))^{3\ell}},\quad (c_{\ell}:=
c_{\ell,0}>0,\,d_{\ell}:=-c_{\ell,1}>0)\quad(\forall\, n\ge 1).
\endmultline\tag 2.2.4
$$
(We write $\sum_ja_jx^j\le_c\sum_jb_jx^j$ when $a_j\le b_j$ for all $j$.) In the same paper he also demonstrated existence of a constant $c>0$ such that 
$$
c_{\ell}\sim c\left(\frac{3}{2}\right)^{\ell}(\ell-1)!,\quad
d_{\ell}\sim c\left(\frac{3}{2}\right)^{\ell}\ell!,\quad (\ell\to\infty).\tag 2.2.5
$$
Later Bagaev and Dmitriev [2] showed that $c=(2\pi)^{-1}$. By now there have been found
other proofs of this fact. See, for instance, Bender et al [3] for an asymptotic
expansion of $c_{\ell}$ due to Meerteens, and \L uczak et al [16] for a rather elementary 
proof based on the behavior of the component size distribution for the critical $G(n,m)$. 
\si

Turn to $E_{r}(x)$, $r\ge 1$. It was shown in [12] that, analogously to (2.2.3),
$$
\aligned
E_{r}(x)=&\sum_{d=0}^{5r}\frac{\e_{r,d}}{(1-T(x))^{3r-d}},\\
\e_{r,d}=&\frac{(6r-2d)!Q_d(r)}{2^{5r}3^{2r-d}(3r-d)!(2r-d)!},
\endaligned\tag 2.2.6
$$
where $Q_0(r)=1$, and, for $d>0$, $Q_d(r)$ is a polynomial of degree $d$.
By Stirling's formula,
$$
\e_{r}:=\e_{r,0}\sim (2\pi)^{-1/2}\left(\frac{3}{2}\right)^{r}r^{r-1/2}e^{-r},
\quad r\to\infty.\tag 2.2.7
$$ 
Formally differentiating both sides of (2.1.8) with respect to $w$ and equating 
coefficients by $w^{\ell-1}$, we get a recurrence relation
$$
r E_{r}(x)=\sum_{k=1}^{r}kC_k(x)E_{r-k}(x). \tag 2.2.8
$$
By (2.2.3) and (2.2.6), the highest power of $(1-T(x))^{-1}$ on both sides of (2.2.8)
is $3r$, and equating the two coefficients we get a recurrence relation involving
$\e_{r}$ and $c_{r}$,
$$
r\, \varepsilon_{r}=\sum_{k=1}^{r}kc_k\e_{r-k},\quad r\ge 1.\tag 2.2.9
$$
\bi

With these preliminaries out of the way, we turn to the formula (2.1.11) for
$\pr(S_n>0,\,\Cal E_n\le L)$. Notice upfront that, for $L=0$---arising when
$\la\to-\infty$---we simply have
$$
\pr(S_n>0,\,\Cal E_n\le 0)=N(n,p)\,[x^n]e^{H(x)},\quad
H(x)=\frac{q}{p}C_{-1}(x)+\frac{1}{2}C_0(x).\tag 2.2.10
$$
The next Lemma provides a counterpart of (2.1.10) and (2.2.10) for $L\in [1,\infty)$.
\proclaim{Lemma 2.2.1} Given $L\in [1,\infty)$,
$$
\pr(S_n>0,\,\Cal E_n\le L)
=N(n,p)\,\sum_{r=0}^{\infty}\left(\frac{p}{2q}\right)^r[x^n]\left\{e^{H(x)}F_r^L(x)\right\},
\tag 2.2.11
$$
where $\{F_r^L(x)\}$ is determined by a recurrence relation
$$
rF_r^L(x)=\frac{1}{2}\sum_{k=1}^{r\wedge L}k\,C_k(x)F_{r-k}^L(x),\quad r\ge 1,\tag 2.2.12
$$
and $F_0^L(x)=1$. (Here $a\wedge b:=\min\{a,b\}$.)
\endproclaim
\si

{\bf Proof of Lemma 2.2.1.\/} Clearly 
$$
\exp\left(\frac{1}{2}\sum_{\ell=1}^L
w^{\ell}C_{\ell}(x)\right)=\sum_{r=0}^{\infty}w^rF_r^L(x), \tag 2.2.13
$$
where $F_r^L(x)$ are some power series, with nonnegative coefficients, convergent for $|x|<e^{-1}$.
This identity implies that
$$
\exp\left(\sum_{\ell=1}^L
w^{\ell}C_{\ell}(x)\right)=\left(\sum_{r=0}^{\infty}w^rF_r^L(x)\right)^2.
$$
Differentiating this with respect to $w$ and replacing $\exp\left(\sum_{\ell=1}^L
w^{\ell}C_{\ell}(x)\right)$ on the left of the resulting identity with 
$\left(\sum_{s=0}^{\infty}w^sF_s^L(x)\right)^2$, we get , after multiplying by $w$,
$$
\left(\sum_{s=0}^{\infty}w^sF_s^L(x)\right)\left(\sum_{\ell=1}^L \ell w^{\ell}C_{\ell}
(x)\right)=2\sum_{r=1}^{\infty}rw^rF_r^L(x).
$$
Equating the coefficients by $w^r$, $r\ge 1$, of the two sides we obtain the recurrence 
(2.2.12).  
\bi

The recurrence (2.2.12) yields a very useful information about $F_r^L(x)$.
\proclaim{Lemma 2.2.2} Let $L>0$. For $r\ge 0$,
$$
F_r^L(x)=\sum_{d=0}^{5r}\frac{f_{r,d}^L}{(1-T(x))^{3r-d}},\tag 2.2.14
$$
and, denoting $f_r^L=f_{r,0}^L$, $g_r^L=-f_{r,1}^L$
$$
\frac{f_r^L}{(1-T(x))^{3r}}-\frac{g_r^L}{(1-T(x))^{3r-1}}\le_c
F_r^L(x)\le_c\frac{f_r^L}{(1-T(x))^{3r}}.\tag 2.2.15
$$
Furthermore the leading coefficients $f_r^L$, $g_r^L$ satisfy a recurrence relation
$$
\align
rf_r^L=&\frac{1}{2}\sum_{k=1}^{r\wedge L}k\, c_k\,f_{r-k}^L;\quad f_0^L=1,\tag 2.2.16\\
r  g_r^L=&\frac{1}{2}\sum_{k=1}^{r\wedge L}k\, c_k\,g_{r-k}^L
+\frac{1}{2}\sum_{k=1}^{r\wedge L}k\, d_k\,f_{r-k}^L;\quad g_0^L=0,\tag 2.2.17
\endalign
$$
so, in particular, $f_r^L>0$ and $g_r^L>0$ for $r>0$.
\endproclaim
\si

{\bf Note.\/} 1. This Lemma and its proof are similar to those  for the generating functions $E_r(x)$ 
obtained in [10]. 
\si

{\bf Proof of Lemma 2.2.2.\/} {\bf (a)\/} We prove (2.2.14) by induction on $r$. (2.2.14) 
holds
for $r=0$ as $F_0^L(x)\equiv 1$ and $f_{0,0}^L=f_0^L=1$. Further, by (2.2.12) and (2.2.3),
$$
F_1^L(x)=\frac{1}{2}C_1(x)=\frac{1}{2}\sum_{d=0}^5\frac{c_{1,d}}{(1-T(x))^{3-d}},
$$
i. e. (2.2.14) holds for $r=1$ too. Assume that $r\ge 2$ and that (2.2.14) holds for
for $r^\prime\in [1,r-1]$. Then, by (2.2.12), (2.2.3) and inductive assumption,
$$
\align
F_r^L(x)=&\frac{1}{2r}\sum_{k=1}^{r\wedge L}kC_k(x)F_{r-k}^L(x)\\
=&\frac{1}{2r}\sum_{k=1}^{r\wedge L}k\sum_{d=0}^{3k+2}
\frac{c_{k,d}}{(1-T(x))^{3k-d}}\sum_{d_1=0}^{5(r-k)}\frac{f_{r-k,d_1}^L}
{(1-T(x))^{3(r-k)-d_1}}\\
=&\frac{1}{2r}\sum_{k=1}^{r\wedge L}k\sum\limits_{d\le 3k+2,\, d_1\le 5(r-k)}
\frac{c_{k,d}\,\,f_{r-k,d_1}^L}{(1-T(x))^{3r-(d+d_1)}}.
\endalign
$$
Here
$$
0\le d+d_1\le 3k+2 +5(r-k)=5r-2(k-1)\le 5r,
$$
so (2.2.14) holds for $r$ as well.
\si

{\bf (b)\/} Plugging (2.2.14) and (2.2.3) into (2.2.12) we get
$$
\multline
\sum_{d=0}^{5r}\frac{f_{r,d}^L}{(1-T(x))^{3r-d}}\\
=\sum_{k=1}^{r\wedge L}\frac{k}{2r}\sum_{d_1=0}^{3k+2}\frac{c_{k,d_1}}
{(1-T(x))^{3k-d_1}}\sum_{d_2=0}^{5(r-k)}\frac{f_{r-k,d_2}^L}
{(1-T(x))^{3(r-k)-d_2}}.
\endmultline
$$
Equating the coefficients by $(1-T(x))^{-3r}$  (by $(1-T(x))^{-3r+1}$ resp.) on the right and 
on the left, we obtain (2.2.16) ((2.2.17) resp.).
\si

{\bf (c)\/} For $r=0$, (2.2.15) holds trivially. For $r\ge 1$, inductively we have: by 
(2.2.4) (upper
bound) and (2.2.12), (2.2.16),
$$
\align
F_r^L(x)\le_c&\,\frac{1}{2r}\sum_{k=1}^{r\wedge L}k\frac{c_k}{(1-T(x))^{3k}}\frac{f_{r-k}^L}
{(1-T(x))^{3(r-k)}}\\
=&\,\frac{1}{(1-T(x))^{3r}}\frac{1}{2r}\sum_{k=1}^{r\wedge L}kc_kf_{r-k}^L\\
=&\,\frac{f_r^L}{(1-T(x))^{3r}};
\endalign
$$
furthermore, by (2.2.4) (lower bound), (2.2.12) and (2.2.16)-(2.2.17),
$$
\align
F_r^L(x)\ge_c&\,\frac{1}{2r}\sum_{k=1}^{r\wedge L}k\left[\frac{c_k}{(1-T(x))^{3k}}-
\frac{d_k}{(1-T(x))^{3k-1}}\right]F_{r-k}^L(x)\\
\ge_c&\,\frac{1}{2r}\sum_{k=1}^{r\wedge L}k\frac{c_k}{(1-T(x))^{3k}}
\left[\frac{f_{r-k}^L}{(1-T(x))^{3(r-k)}}-\frac{g_{r-k}^L}{(1-T(x))^{3(r-k)-1}}
\right]\\
&\,-\frac{1}{2r}\sum_{k=1}^{r\wedge L}k\frac{d_k}{(1-T(x))^{3k-1}}
\cdot\frac{f_{r-k}^L}
{(1-T(x))^{3(r-k)}}\\
=&\,\frac{f_r^L}{(1-T(x))^{3r}}-\frac{1}{(1-T(x))^{3r-1}}\left[\frac{1}{2r}
\sum_{k=1}^{r\wedge L}k c_kg_{r-k}^L+\frac{1}{2r}
\sum_{k=1}^{r\wedge L}k d_k f_{r-k}^L\right]\\
=&\,\frac{f_r^l}{(1-T(x))^{3r}}-\frac{g_r^L}{(1-T(x))^{3r-1}}.
\endalign
$$
\qed
\bi

To make the bound (2.2.15) work we need to have a close look at the sequence
$\{f_r^L,g_r^L\}_{r\ge 0}$. First of all, it follows from (2.2.16) that
$$
f_r^{L}\le f_r:=f_r^{\infty},\quad g_r^L\le g_r:=g_r^{\infty}.
$$
That is $f_r$ and $-g_r$ are the coefficients by $(1-T(x))^{-3r}$ and $(1-T(x))^
{-3r+1}$ in the expansion (2.2.13) for $F_r(x):=F_r^{\infty}(x)$.
 Now, using (2.2.13) for $L=\infty$ and (2.1.8),
we see that
$$
\left(\sum_{r\ge 0}w^rF_r(x)\right)^2=\sum_{r\ge 0}w^r E_r(x).
$$
So, equating the coefficients by $w^r$, $r\ge 0$, we get 
$$
\sum_{k=0}^rF_k(x)F_{r-k}(x)=E_r(x).
$$
Plugging (2.2.6) and (2.2.14) (with $L=\infty$), and comparing coefficients by $(1-T(x))^{-3r}$
($(1-T(x))^{-3r+1}$, resp.), we obtain  
$$
\sum_{k=0}^rf_kf_{r-k}=\e_{r,0};\quad 2\sum_{k=0}^rf_kg_{r-k}=-\e_{r,1}.
$$ 
In particular,
$$
f_r\le \frac{1}{2}\e_{r,0},\quad g_r\le -\frac{1}{2}\e_{r,1}.
$$
Consequently, using (2.2.6) for $r\ge 2$ and $d=0$,
$$
\align
f_r=&\frac{1}{2}\varepsilon_{r,0}-\frac{1}{2}\sum_{k=1}^{r-1}f_kf_{r-k}
\ge\frac{1}{2}\varepsilon_{r,0}-\frac{1}{2}\sum_{k=1}^{r-1}\frac{1}{2}\varepsilon_{k,0}
\frac{1}{2}\varepsilon_{r-k,0}\\
\ge& \frac{\e_{r,0}}{2}\left(1-\frac{1}{4}\sum_{j=1}^{r-1}\binom{r}{j}^{-1}
\frac{\binom{r}{j}\binom{2r}{2j}\binom{3r}{3j}}{\binom{6r}{6j}}\right)\\
\ge&\frac{\e_{r,0}}{2}\left(1-\frac{1}{4}\sum_{j=1}^{r-1}\binom{r}{j}^{-1}\right)\\
\ge&\frac{\e_{r,0}}{2}(1-1/r),
\endalign
$$
that is
$$
\frac{\e_{r,0}}{2}(1-1/r)\le f_r\le \frac{\e_{r,0}}{2}\sim\frac{1}
{2\sqrt{2\pi}}\left(\frac{3}{2}\right)^rr^{r-1/2}e^{-r}, \, (r\to\infty),\tag 2.2.18
$$
see (2.2.7).
Furthermore, using (2.2.6) for $r>0$ and $d=1$,
$$
g_r\le_b\left(\frac{3}{2}\right)^r r^{r+1/2}e^{-r}.\tag 2.2.19
$$
And one can prove a matching lower bound for $g_r$. Hence, like $\varepsilon_r$, $f_r$, $g_r$ grow essentially as  $r^r$,
too fast for $F_r(x)=F_r^{\infty}(x)$ to be useful for asymptotic estimates. The next Lemma
(last in this subsection) shows that, in a pleasing contrast, $f_r^L$, $g_r^L$ grow much 
slower when
$r\gg L$.
\proclaim{Lemma 2.2.3} There exists $L_0$  such that, for $L\ge L_0$, 
$$
f_r^L\le_b\left(\frac{3L}{2e}\right)^r,\quad g_r^L\le_br\left(\frac{3L}{2e}\right)^r,\quad\forall\,r\ge 0.\tag 
2.2.20
$$
\endproclaim
\si

{\bf Proof of Lemma 2.2.3.\/} {\bf (a)\/}  It is immediate from (2.2.18),  (2.2.19)
that, for some absolute constant $A$ and all $L>0$, 
$$
f_r^L=f_r\le A\left(\frac{3L}{2e}\right)^r,\quad g_r^L=g_r\le Ar\left(\frac{3L}{2e}\right)^r\quad 0\le r \le L.
$$
Let us prove existence  an integer $L>0$, with a property: if for some $s\ge L$ and all $t\le s$,
$$
f_t^L\le A\left(\frac{3L}{2e}\right)^t,\quad g_t^L\le At\left(\frac{3L}{2e}\right)^t, \tag 2.2.21
$$
then
$$
f_{s+1}^L\le A\left(\frac{3L}{2e}\right)^{s+1},\quad g_{s+1}^L\le A(s+1)\left(\frac{3L}{2e}\right)^{s+1}.
$$
By (2.2.16),  (2.2.21), and (2.2.5), there exists an absolute constant $B>0$ such that
$$
(s+1)f_{s+1}^L\le AB\left(\frac{3L}{2e}\right)^{s+1}L^{1/2}\sum_{k=1}^L\left(\frac{k}{L}\right)^k.
$$
A function $(x/L)^x$ attains its minimum on $[0,L]$ at $x=L/e$, and it is easy to show
that
$$
(x/L)^x\le \left\{\alignedat2
&e^{-x},\quad&&  x\le L/e,\\
&e^{-(L-x)(3-e)/2},\,\,\quad&& x\ge L/e.\endalignedat\right.
$$
Since $s+1\ge L$, we obtain then
$$
\align
f_{s+1}^L\le&\, AB\left(\frac{1}{1-e^{-1}}+\frac{1}{1-e^{-(3-e)/2}}\right)\cdot L^{-1/2}\left(\frac{3L}{2e}
\right)^{s+1}\\
\le&\, A\left(\frac{3L}{2e}\right)^{s+1},
\endalign
$$
if we choose
$$
L\ge L_1:=B^2\left(\frac{1}{1-e^{-1}}+\frac{1}{1-e^{-(3-e)/2}}\right)^2.
$$
Likewise, by (2.2.17), (2.2.21) and (2.2.5),
$$
\align
(s+1)g_{s+1}^L\le& AB (s+1)\left(\frac{3L}{2e}\right)^{s+1}L^{1/2}\sum_{k=1}^L\left(\frac{k}{L}\right)^k\\
&+AB^\prime (s+1)\left(\frac{3L}{2e}\right)^{s+1}L^{1/2}\sum_{k=1}^L\left(\frac{k}{L}\right)^k,
\endalign
$$
so that 
$$
g_{s+1}^L\le A(s+1)\left(\frac{3L}{2e}\right)^{s+1},
$$
if we choose
$$
L\ge L_2:=(B+B^\prime)^2\left(\frac{1}{1-e^{-1}}+\frac{1}{1-e^{-(3-e)/2}}\right)^2.
$$
Thus, picking $L=\max\{L_1,L_2\}=L_2$,  we can accomplish the inductive step, from $s\, (\ge L)$ to $s+1$,
showing that, for this $L$,  (2.2.20) holds for all $t$.
\qed
\bi

Combining (2.2.10), Lemma 2.2.1, Lemma 2.2.2, 
we  bound $\pr(S_n>0,\,\Cal E_n\le L)$.
\proclaim{Proposition 2.2.4} Let $L\in [0,\infty)$. Then
$$
\Sigma_1\le \pr(S_n>0,\,\Cal E_n\le L)\le \Sigma_2.
$$
Here
$$
\aligned
\Sigma_1=&\,N(n,p)\sum_{r\ge 0}\left(\frac{p}{2q}\right)^r[x^n]\,
\left[\frac{f_r^Le^{H(x)}}{(1-T(x))^{3r}}-\frac{g_r^Le^{H(x)}}{(1-T(x))^{3r-1}}\right],\\
\Sigma_2=&\,N(n,p)\sum_{r\ge 0}\left(\frac{p}{2q}\right)^r[x^n]\,
\frac{f_r^Le^{H(x)}}{(1-T(x))^{3r}},
\endaligned\tag 2.2.22
$$
and
$$
f_r^L\left\{\alignedat2
&=f_r,\quad&&r\le L,\\
&\le_b \left(\frac{3L}{2e}\right)^r,\quad&&r\ge L,\endalignedat\right.
\quad
g_r^L\left\{\alignedat2
&=g_r,\quad&&r\le L,\\
&\le_b r\left(\frac{3L}{2e}\right)^r,\quad&&r\ge L,\endalignedat\right.\tag 2.2.23
$$
with $f_r$, $g_r$ satisfying the conditions (2.2.18)-(2.2.19).
\endproclaim
\si

{\bf Note.\/} The relations (2.2.22)-(2.2.23) indeed cover the case $L=0$ since in this case
$f_0=1$, $g_0=0$ and $f_r^L=g_r^L=0$ for $r>0$.
\bi
{\bf 2.3 Asymptotic formula for $\pr(S_n>0)$.\/}
\bi
The Proposition 2.2.4 makes it clear that we need to find an asymptotic formula for 
$$
N(n,p)\phi_{n,w},\quad \phi_{n,w}:=[x^n]\,\frac{e^{H(x)}}{(1-T(x))^{w}},\quad w=0,3,6\dots \tag 2.3.1
$$
Using $N(n,p)=n!q^{n^2/2}(pq^{-3/2})^n$ and Stirling's formula for $n!$, with some
work we obtain
$$
\multline
N(n,p)=\sqrt{2\pi n}\exp\biggl[-n\frac{3}{2}+n^{2/3}\frac{\la}{2}-n^{1/3}\frac{\la^2}{2}\\
+\frac{\la^3}{3}+\frac{5}{4}+O(n^{-1/3}(1+\la^4))\biggr].
\endmultline\tag 2.3.2
$$
The big-Oh term here is $o(1)$ if $|\la|=o(n^{1/12})$, which is the condition of Theorem 1.1.
\si

Turn to $\phi_{n,w}$.
Since the function in question is analytic for $|x|<e^{-1}$,
$$
\phi_{n,w}=\frac{1}{2\pi i}\oint_{\Gamma}\frac{e^{H(x)}}{x^{n+1}(1-T(x))^{w}}\,dx,
$$
where $\Gamma$ is a simple closed contour enclosing the origin and lying in the
disc $|x|<e^{-1}$. By (2.1.10),  (2.2.1)-(2.2.2), the function in (2.3.1) depends on $x$ only
through $T(x)$, which satisfies $T(x)=xe^{T(x)}$. This suggests introducing a new variable
of integration $y$, such that $ye^{-y}=x$, i. e. 
$$
y=T(x)=\sum_{k\ge 1}\frac{x^k}{k!}k^{k-1},\quad |x|<e^{-1}.
$$
Picking a simple closed contour $\Gamma^\prime$ in the $y$-plane such that its image
under $x=ye^{-y}$ is a simple closed contour $\Gamma$ within the disc $|x|<e^{-1}$, and using
(2.2.1)-(2.2.2), we obtain
$$
\aligned
\phi_{n,w}=&\frac{1}{2\pi i}\oint_{\Gamma^\prime}y^{-n-1}e^{ny}\exp\left(\kappa (y)-\frac{y}{4}-
\frac{y^2}{8}\right)
(1-y)^{3/4-w}\,dy,\\
\kappa(y):=&\frac{q}{p}\left(y-\frac{y^2}{2}\right);
\endaligned\tag 2.3.3
$$
$q/p\sim n$, so $y^{-n}e^{ny}e^{\kappa(y)}$ would have fully accounted  for asymptotic 
behavior of the integral, had it not been  for the factor $(1-y)^{3/4-w}$.  
Once $\Gamma^\prime$ is picked, it can be replaced by any
circular contour $y=\rho e^{i\theta}$, $\theta\in (-\pi,\pi]$, $\rho<1$. (The condition
$\rho<1$ is dictated by the factor $(1-y)^{3/4-w}$.) And (2.3.3)
becomes
$$
\aligned
\phi_{n,w}=&\frac{1}{2\pi}I(w),\\
I(w):=&\int_{-\pi}^{\pi}e^{h(\rho,\theta)}\exp\left(-\rho e^{i\theta}/4-\rho^2e^{i2\theta}/8
\right)(1-\rho e^{i\theta})^{3/4-w}\,d\theta,\\
h(\rho,\theta)=&\frac{q}{p}\left(\rho e^{i\theta}-\frac{\rho^2e^{i2\theta}}{2}\right)
+n\rho e^{i\theta}-n(\ln\rho +i\theta).
\endaligned\tag 2.3.4
$$
Let us choose $\rho<1$ in such a way that, as a function of $\theta$, $|e^{h(\rho,\theta)}|$ 
attains its maximum at $\theta=0$. Now $|e^{h(\rho,\theta)}|=e^{f(\rho,\theta)}$, with
$$
f(\rho,\theta)=\text{Re }h(\rho,\theta)=\frac{q}{p}\,\rho\cos\theta-\frac{q}{2p}\,\rho^2\cos 2\theta+n\rho
\cos\theta-n\ln\rho,
$$
so that
$$
f^\prime_{\theta}(\rho,\theta)=\frac{2q}{p}\,\rho^2\sin\theta\left(\cos\theta-\frac{1+np/q}{2\rho}\right).
$$
Then $f^\prime_{\theta}(\rho,\theta)>0$ ($<0$ resp.) for $\theta<0$ ($\theta>0$ resp.) if
$$
\rho<\frac{1}{2}(1+np/q).\tag 2.3.5
$$
Let us set $\rho=e^{-an^{-1/3}}$, where $a=o(n^{1/3})$, since we want $\rho
\to 1$. Now
$$
\frac{1}{2}(1+np/q)>1+\frac{\la}{2}n^{-1/3},\quad\rho\le 1-an^{-1/3}+\frac{a^2}{2}
n^{-2/3};
$$
so (2.3.5) is obviously satisfied if  
$$
a+\frac{\la}{2}\ge \frac{a}{2}.\tag 2.3.6
$$
(2.3.6) is trivially met if $\la\ge 0$. For $\la<0$, $|\la|=o(n^{1/3})$, 
(2.3.6) is met if $a\ge |\la|$. In all cases we will assume that
$\liminf a>0$.
\si

Why do we want $a=o(n^{1/3})$? Because, as a function of $\rho$, $h(\rho,0)$ attains its
minimum at $np/q\sim 1$, if $\la<0$ is fixed, and in this case $np/q<1$, and the minimum point is 
$1$ if $\la\ge 0$. So our $\rho$ is a reasonable approximation of the saddle point of
$|h(\rho,\theta)|$, dependent on $\la$,
chosen from among the {\it feasible\/} values, i. e. those strictly below $1$. Characteristically
$\rho$ is very close to $1$, the singular point of the factor $(1-y)^{3/4-w}$, which is especially
influential for large $w$'s. Its presence rules out a ``pain-free'' application of general tools such as
Watson's Lemma.
\si

Under (2.3.6), 
$$
|f^\prime_{\theta}(\rho,\theta)|\ge \frac{a}{2}n^{2/3}|\sin\theta|,
$$
and $\text{sign}f^\prime_{\theta}(\rho,\theta)=-\text{sign }\theta$, so that
$$
\aligned
f(\rho,\theta)\le&\, f(\rho,0)-\frac{a}{2}n^{2/3}\int_0^{|\theta|}\sin z\,dz\\
=&\,f(\rho,0)-an^{2/3}\sin^2(\theta/2)\\
\le&\,f(\rho,0)-a\pi^{-2}n^{2/3}\theta^2=h(\rho,0)-a(\pi)^2n^{2/3}\theta^2.
\endaligned\tag 2.3.7
$$
Let us break the integral $I(w)$ in (2.3.4) into two parts, $I_1(w)$ for $|\theta|\le
\theta_0$, and $I_2(w)$ for $|\theta|\ge \theta_0$, where
$$
\theta_0=\pi n^{-1/3}\ln n.
$$
Since $f(\rho,\theta)$ is decreasing with $|\theta|$, and $|1-\rho e^{i\theta}|\ge 1-\rho$,
it follows from  (2.3.7) that
$$
\aligned
|I_2(w)|\le_b&\, \bigl(1-e^{-an^{-1/3}}\bigr)^{-w}e^{f(\rho,\theta_0)}\\
\le&\,\bigl(a_1n^{-1/3}\bigr)^{-w}e^{h(\rho,0)}\exp(-a\ln^2 n);\\
a_1:=&\,n^{1/3}\bigl(1-e^{-an^{-1/3}}\bigr).
\endaligned\tag 2.3.8
$$
Turn to $I_1(w)$. This time $|\theta|\le \theta_0$. First, let us  write
$$
\rho e^{i\theta}=e^{-sn^{-1/3}},\quad s=a-it,\quad t:=n^{1/3}\theta;
$$
so $|s|\le a+\pi \ln n$. The second (easy) exponent in the integrand of $I_1(w)$ is 
asymptotic to $-3/8$, or more precisely,
$$
\aligned
-\frac{1}{4}e^{-sn^{-1/3}}-\frac{1}{8}e^{-2sn^{-1/3}}=&Q_2(a)+O(|t|n^{-1/3}),\\
Q_2(a):=-\frac{1}{4}e^{-an^{-1/3}}-\frac{1}{8}e^{-2an^{-1/3}}.
\endaligned\tag 2.3.9
$$
Determination of a usable
asymptotic formula for $h(\rho,\theta)$ is more laborious. It is convenient to
set  $q/p=n e^{-\mu n^{-1/3}}$; thus
$$
\mu=n^{1/3}\ln\frac{np}{q}\ge n^{1/3}\ln (1+\la n^{-1/3})\ge \la(1-\la n^{-1/3}/2),
$$
and 
$$
\mu-\la=O(n^{-2/3}+n^{-1/3}\la^2).
$$
Using the new parameters $s$ and $\mu$
we transform  the formula (2.3.4) for $h(\rho,\theta)$ to
$$
h(\rho,\theta)= 
n\left(e^{-(\mu+s)n^{-1/3}}-\frac{1}{2}e^{-(\mu+2s)n^{-1/3}}+e^{-sn^{-1/3}}+sn^{-1/3}
\right).
$$
Approximating the three exponents by
the $4$-th degree Taylor polynomials, we obtain
$$
\aligned
h(\rho,\theta)
=&\,n\left[\frac{3}{2}-n^{-1/3}\,\frac{\mu}{2} +n^{-2/3}\,\frac{\mu^2}{4}
-n^{-1}\frac{\mu^3}{12}\right]+Q_1(\mu,a)\\
&+\left(\frac{\mu s^2}{2}+\frac{s^3}{3}\right)+O\bigl(D_1(t)\bigr);\\
Q_1(\mu,a):=&n^{-1/3}\left[\frac{(\mu+a)^4}{4!}-\frac{(\mu+2a)^4}{4!2}+
\frac{a^4}{4!}\right];\\
D_1(t):=&n^{-1/3}|t|(|\la|+a+\ln n)^3+n^{-2/3}(|\la|+a+\ln n)^5.
\endaligned\tag 2.3.10
$$
(Explanation: the second summand in $D_1(t)$ is 
the approximation error bound for
each of the Taylor polynomials; the first summand is the common bound of 
$|(\mu+a)^4-(\mu+s)^4|$,
$|(\mu+2s)^4-(\mu+2a)^4|$, and $|s^4-a^4|$, times $n^{-1/3}$.)  And we notice 
immediately 
that both $Q_1(\mu,a)$
and $D_1(t)$ are $o(1)$ if, in addition to $|\la|=o(n^{1/12})$, we require
that $a =o(n^{1/12})$ as well, a condition we assume from now on. Obviously 
$O(D_1(t))$ absorbs the remainder term $O(|t|n^{-1/3})$ from (2.3.9). 
\si

Furthermore, since 
$$
n^{-1/3}\mu=\ln\left(\frac{np}{q}\right)=n^{-1/3}\la-n^{-2/3}\frac{\la^2}{2}+n^{-1}
\left(\frac{\la^3}{3}+1\right)+O(n^{-4/3}(1+\la^4)),
$$
for the cubic polynomial of $n^{-1/3}\mu$ in (2.3.10) we have
$$
\align
&n\left[\frac{3}{2}-n^{-1/3}\,\frac{\mu}{2} +n^{-2/3}\,\frac{\mu^2}{4}
-n^{-1}\frac{\mu^3}{12}\right]\\
=&\,n\frac{3}{2}-n^{2/3}\frac{\la}{2}+n^{1/3}\frac{\la^2}{2}-\frac{\la^3}{2}-
\frac{1}{2}+O(n^{-1/3}(1+\la^4)).
\endalign
$$
Observe that the first three summands are those in the exponent of the
formula (2.3.2) for $N(n,p)$ times $(-1)$.) 
Therefore, using (2.3.2) for $N(n,p)$,
$$
\aligned
&N(n,p)\exp\left(h(\rho,\theta)-\frac{1}{4}\rho e^{i\theta}-\frac{1}{8}\rho^2e^{i2\theta}
\right)\\
=&\bigl[1+O(D_1(t))\bigr]\sqrt{2\pi n}
\cdot\exp\left(-\frac{\la^3}{6}+\frac{3}{4}+Q(\mu,a)+
\frac{\mu s^2}{2}+\frac{s^3}{3}\right);\\
Q(\mu,a):=&Q_1(\mu,a)+Q_2(a)+O(n^{-1/3}(1+\la^4)),
\endaligned\tag 2.3.11
$$
and $Q(\mu,a)=o(1)$ as $\la,\,a=o(n^{1/12})$. In particular, using (2.3.8), (2.3.10) for 
$\theta=0$, i. e. $s=a$, we see that
$$
N(n,p)|I_2(w)|\le_bn^{1/2} (a_1n^{-1/3})^{-w}e^{-a\ln^2 n}\exp\left(-
\frac{\la^3}{6}+\frac{\mu a^2}{2}+\frac{a^3}{3}\right).\tag 2.3.12
$$
Furthermore,  switching integration from $\theta$ to $t=n^{1/3}\theta$, 
the contribution of the remainder term $O(D_1(t))$ to $N(n,p)I_1(w)$ is 
$O(\delta_{n,w})$,
$$
\delta_{n,w}:=n^{-1/12}\,\frac{(a+\ln n)^{3/4}e^{-\la^3/6}}{(a_1n^{-1/3})^{w}}
\int_{-\infty}^{\infty}
\left|\exp\left(\frac{\mu s^2}{2}+\frac{s^3}{3}\right)\right|\,D_1(t)\,dt.
$$
(Explanation: $n^{-1/12}=n^{1/2}n^{-1/3} n^{-1/4}$, with $n^{-1/4}$ coming from 
$n^{-1/4}(a+\pi\ln n)^{3/4}$, an upper bound of $|1-\rho e^{i\theta}|^{3/4}$, for
$|\theta|\le \theta_0$.) 
\si

Now
$$
\left|\exp\left(\frac{\mu s^2}{2}+\frac{s^3}{3}\right)\right|
=\exp\left[\frac{\mu a^2}{2}+\frac{a^3}{3}-\left(\frac{\mu}{2}+a\right)t^2\right],
$$
where, see (2.3.6), 
$$
\frac{\mu}{2}+a=\frac{\la}{2}+a+O(n^{-2/3}+n^{-1/3}\la^2)>0,
$$
since $\liminf a>0$, and $a\ge |\la|$ if $\la<0$. Hence, see (2.3.10) for $D_1(t)$, we have 
$\delta_{n,w}\le_b\Delta_{n,w}$, where
$$
\multline
\Delta_{n,w}:=\,
n^{-1/12+w/3}\cdot a_1^{-w}\exp\left(-\frac{\la^3}{6}+\frac{\mu a^2}{2}+\frac{a^3}{3}
\right)\\
\cdot (a+\ln n)^{3/4} \left[n^{-1/3}\,\frac{(|\la|+a+\ln n)^3}{\mu/2+a}
+n^{-2/3}\,\frac{(|\la|+a+\ln n)^5}{(\mu/2+a)^{1/2}}\right].
\endmultline\tag 2.3.13
$$
The denominators $\mu/2+a$, $(\mu/2+a)^{1/2}$  come from the integrals
$$
\int_{-\infty}^{\infty}|t|^k\exp\left[-\left(\frac{\mu}{2}+a\right)t^2
\right]\,dt =c_k(\mu/2+a)^{-(k+1)/2},\quad (k\ge 0),
$$
for $k=0,1$. Clearly $\Delta_{n,w}$ absorbs the bound (2.3.12). 
\si

Thus, switching from $\theta$ to $s=a-in^{1/3}\theta$, it remains 
to evaluate sharply 
$$
\aligned
&-i(2\pi)^{-1/2}n^{1/6}\exp\left(-\frac{\la^3}{6}+\frac{3}{4}+Q(\mu,a)\right)\\
&\cdot
\int\limits_{s_1}
^{s_2}
\exp\left(\frac{\mu s^2}{2}+\frac{s^3}{3}\right)(1-e^{-sn^{-1/3}})^{3/4-w}\,ds;
\endaligned\tag 2.3.14
$$
here $s_1=a-in^{1/3}\theta_0$, $s_2=a+in^{1/3}\theta_0$, and the integral is over the
vertical line segment connecting $s_1$ and $s_2$. Lastly we need to 
estimate an error coming from replacing  $(1-e^{-sn^{-1/3}})^{3/4-w}$ with
a genuinely palatable $(sn^{-1/3})^{3/4-w}$.  Using
$$
\align
&|sn^{-1/3}|\ge |1-e^{-sn^{-1/3}}|\ge |1-e^{-an^{-1/3}}|,\quad (s=a-it),\\
&|x^u-1|\le u|x-1|,\quad (u\ge 1,\,|x|\le 1),
\endalign
$$
we have: for $u\ge 1$
$$
\align
&\left|\frac{1}{(1-e^{-s n^{-1/3}})^u}-\frac{1}{(s n^{-1/3})^u}\right|
\le \frac{1}{|1-e^{-a n^{-1/3}}|^u}\left|1-\left(\frac{1-e^{-s n^{-1/3}}}
{sn^{-1/3}}\right)^u\right|\\
\le& \frac{u}{|1-e^{-a n^{-1/3}}|^u}\left|1-\frac{1-e^{-s n^{-1/3}}}
{sn^{-1/3}}\right|\le_b \frac{u\,|sn^{-1/3}|}{|1-e^{-a n^{-1/3}}|^u}
\le \frac{u(a+|t|)n^{-1/3}}{|1-e^{-a n^{-1/3}}|^u}.
\endalign
$$
Also, for $s$ in question,
$$
|1-e^{-sn^{-1/3}}|\ge 0.5 |sn^{-1/3}|.
$$
So
$$
\align
|(1-e^{-sn^{-1/3}})^{3/4}-(sn^{-1/3})^{3/4}|=&|sn^{-1/3}|^{3/4}
\left|\left(\frac{1-e^{-sn^{-1/3}}}{sn^{-1/3}}\right)^{3/4}
-1\right|\\
\le_b&|sn^{-1/3}|^{3/4+1}\le_bn^{-7/12}(a+|t|)^{7/4}.
\endalign
$$
Combining these two estimates, we have: for $w\in \{0,3,4,\dots\}$,
$$
|(1-e^{-sn^{-1/3}})^{3/4-w}-(sn^{-1/3})^{3/4-w}|\le_b
(w+1)\frac{n^{-7/12}(a+|t|)^{7/4}}{(a_1n^{-1/3})^w};
$$
see (2.3.8) for $a_1$. Consequently, replacing $(1-e^{-sn^{-1/3}})^{3/4-w}$
in (2.3.14) with $(sn^{-1/3})^{3/4-w}$ incurs an
{\it additive\/} error of order 
$$
(w+1)n^{-1/12+w/3}\cdot a_1^{-w}\exp\left(-\frac{\la^3}{6}+\frac{\mu a^2}
{2}+\frac{a^3}{3}\right)\cdot n^{-1/3}a^{5/4},
$$
at most; thus the error is easily $O((w+1)\Delta_{n,w})$, see (2.3.13) for 
$\Delta_{n,w}$.
\bi

While these bounds will suffice for $\la=O(1)$, the case $\la\to\infty$
requires a sharper approximation of $(1-e^{-sn^{-1/3}})^{3/4-w}$  for $w=O(\la^3)$. We write
$$
\aligned
(1-e^{-sn^{-1/3}})^{3/4-w}=&\,(sn^{-1/3})^{3/4-w}\exp\left[(3/4-w)
\ln\frac{1-e^{-sn^{-1/3}}}{sn^{-1/3}}\right]\\
=&\,(sn^{-1/3})^{3/4-w}\exp\bigl[Q_3(w,a)+O(D_3(w,t))\bigr];\\
Q_3(w,a):=&\,(3/4-w)
\ln\frac{1-e^{-an^{-1/3}}}{an^{-1/3}};\\
D_3(w,t):=&\,(w+1)tn^{-1/3}.
\endaligned\tag 2.3.15
$$
Notice that $Q_3(a,w)\to 0$ as $wa=O(\la^3 n^{1/12})=o(n^{1/3})$, and $D_3(t,w)\to 0$
as $w\ln n=o(n^{1/3})$.  The expression
(2.3.14) therefore becomes
$$
\aligned
-i(2\pi)^{-1/2}&n^{-1/12+w/3}\exp\left(-\frac{\la^3}{6}+\frac{3}{4}+Q(\mu,w,a)\right)\\
&\cdot
\int\limits_{s_1}
^{s_2}
\exp\left(\frac{\mu s^2}{2}+\frac{s^3}{3}\right)s^{3/4-w}\,ds+O((w+1)
\tilde\Delta_{n,w});\\
Q(\mu,w,a):=&Q(\mu,a)+Q_3(w,a);\\
\tilde \Delta_{n,w}:=&n^{-5/12+w/3}\cdot a^{-w}\exp\left(-\frac{\la^3}{6}
+\frac{\mu a^2}{2}+\frac{a^3}{3}\right)(a+\ln n)^{3/4}.
\endaligned\tag 2.3.16
$$
\si

Finally, after this replacement we can
extend the integration to $(a-i\infty,a+i\infty)$, since the attendant additive error 
is easily shown to be absorbed by $(w+1)\Delta_{n,w}$ for all $w$, and
by $(w+1)\tilde\Delta_{n,w}$ if $w=O(\la^3)$.  
\proclaim{Lemma 2.3.1} Suppose that $\lambda=o(n^{1/12})$. Let $a\ge |\la|$ be such 
that $\lim a>0$, $a=o(n^{1/12})$. Then, denoting $\mu=n^{1/3}\ln(np/q)$, 
$$
\multline
N(n,p)\,[x^n]\frac{e^{H(x)}}{(1-T(x))^w}\\
=-i(2\pi)^{-1/2}e^{3/8+o(1)}n^{-1/12+w/3}e^{-\mu^3/6}\!\!\int\limits_{a-i\infty}^
{a+i\infty}\!\!\!s^{3/4-w}\exp\!\!\left(\frac{\mu s^2}{2}+\frac{s^3}{3}\right)ds\\
+O((w+1)R_{n,w}),
\endmultline\tag 2.3.17
$$
with $R_{n,w}\le\Delta_{n,w}$ for all $w$, and $R_{n,w}=\Delta_{n,w}\wedge
\tilde\Delta_{n,w}$ if $wa$ and $w\ln n$ are both $o(n^{-1/3})$. 
Furthermore, shifting the integration line to
$\{s=b+it:\,t\in (-\infty,\infty)\}$ does not change the value of the integral as 
long as $b\wedge(\mu/2+b)$ 
remains positive.
\endproclaim
\si

{\bf Proof of Lemma 2.3.1.\/} We only have to explain preservation of the integral,
and why $e^{-\la^3/6}$ can be replaced with $e^{-\mu^3/6}$.  
Given such a $b$, pick $T>0$ and 
introduce two horizontal line segments, $C_{1,2}=\{s=\a\pm i T:\,\a\in [a,b]\}$, the top segment
and the bottom segment being respectively right and left oriented. On $C_1\cup C_2$,
$$
\text{Re }\left(\frac{\mu s^2}{2}+\frac{s^3}{3}\right)=\frac{\mu \a^2}{2}+
\frac{\a^3}{3}-T^2\left(\frac{\mu}{2}+\a\right),
$$
and
$$
\frac{\mu}{2}+\a\ge \frac{\mu}{2}+(a\wedge b)>0.
$$
Therefore
$$
\lim_{T\to\infty}\int\limits_{C_1\cup C_2}s^{3/4-w}\exp\left(\frac{\mu s^2}{2}+
\frac{s^3}{3}\right)\,ds=0.
$$
As for $e^{-\la^3/6}\sim e^{-\mu^3/6}$, this follows from
$$
|\la^3-\mu^3|\le_b\la^2(n^{-2/3}+n^{-1/3}\la^2)=(n^{-1/3}\la)^2+n^{-1/3}\la^4\to 0.
\tag 2.3.18
$$
\qed
\bi

In the context of the critical random graph $G(n,m)$, the integral appearing 
in (2.3.17) was encountered and studied in [12]. Following [12], introduce
$$
A(y,\mu)=\frac{e^{-\mu^3/6}}{2\pi i}\int\limits_{a-i\infty}^{a+i\infty}s^{1-y}
\exp\left(\frac{\mu
 s^2}{2}+\frac{s^3}{3}\right)\,ds.\tag 2.3.19
$$
We know that this integral is well defined, and does not depend on $a$, if $a>0$ and 
$a>-\mu/2$. It
was shown in [12] that (1)
$$
A(y,\mu)=\frac{e^{-\mu^3/6}}{3^{(y+1)/3}}\sum_{k=0}^{\infty}
\frac{(3^{2/3}\mu/2)^k}{k!\Gamma((y+1-2k)/3)},\tag 2.3.20
$$
(2) $A(y,\mu)\ge 0$ for $y>0$, $A(y,\mu)>0$ for $y\ge 2$, and (3)
$$
A(y,\mu)\sim \left\{\alignedat2
&(2\pi)^{-1/2}|\mu|^{1/2-y},\quad&&\mu\to-\infty,\\
&\frac{e^{-\mu^3/6}}{2^{y/2}\Gamma(y/2)\mu^{1-y/2}},\quad&&\mu\to\infty.
\endalignedat\right.\tag 2.3.21
$$
We will also need two bounds
$$
\aligned
A(y,\mu)\le_b&\, e^{2|\mu|^3/3}\,\frac{(2/3)^{\frac{y+1}{3}}}{\Gamma\left(\frac{y+1}
{3}\right)},\\
A(y,\mu)\le_b&\,(a+\mu/2)^{-1/2} a^{1-y}\exp\left(-\frac{\mu^3}{6}+\frac{\mu a^2}{2}+
\frac{a^3}{3}\right),
\endaligned\tag 2.3.22
$$
(the second bound holding for $y\ge 1$ and $a+\mu/2>0$), and an asymptotic formula: 
if  $\mu\to\infty$, $y\to\infty$, and $y=O(\mu^3)$, then
$$
A(y,\mu)\sim (2\pi)^{-1/2}(y\xi^{-2}+\mu+2\xi)^{-1/2}\xi^{1-y}\exp\left(-
\frac{\mu^3}{6}+\frac{\mu\xi^2}{2}+
\frac{\xi^3}{3}\right), \tag 2.3.23
$$
where $\xi=\xi(y,\mu)$ is a unique positive root of
$$
\mu\xi^2+\xi^3=y. 
$$
Also, if $y=O(\la^3)$, then
$$
A(y,\mu)\le_b\mu^{-1/2}\xi^{1-y}\exp\left(-\frac{\mu^3}{6}+\frac{\mu\xi^2}{2}+
\frac{\xi^3}{3}\right), \tag 2.3.24
$$
(See Appendix for a proof of (2.3.22) and (2.3.23)-(2.3.24).)
\bi

With  $A(y,\mu)$, we write (2.3.17) more compactly:
$$
\multline
N(n,p)\,[x^n]\frac{e^{H(x)}}{(1-T(x))^w}\\=(2\pi)^{1/2}e^{3/8+o(1)}A(1/4+w,\mu)
n^{-1/12+w/3}
+O((w+1)\Delta_{n,w}).
\endmultline\tag 2.3.25
$$
Let us use (2.3.25) for asymptotic evaluation of $\pr(S_n>0,\,\Cal E_n\le L)$ given 
by (2.2.10)-(2.2.11). 
\bi

{\it Case $|\la|=O(1)$. \/} According to (2.1.16), we can pick $L\to\infty$ as 
slowly as we wish. We pick $L=\ln^{1/4}n$. 
\si

As a first step, let us estimate the overall contributions, $R_n^{(1)}$ and $R_n^{(2)}$, of 
the remainders $O((w+1)R_{n,w})$ to the bounds $\Sigma_1$ and $\Sigma_2$
in Proposition 2.2.3. In this case we choose $a=(L)^{1/3}$ for each
$w$, and $R_{n,w}=\Delta_{n,w}$. Consider $R_n^{(2)}$ first.  By (2.2.19) and (2.3.13), and
dropping $(3r+1)(np/2q)^r=(3r+1)(1/2+o(1))^r$ factor,
$$
\aligned
R_n^{(2)}\le_b&\,n^{-1/12}\cdot\bigl(n^{-1/3}\ln^{15/4}n\bigr)\\
&\cdot \exp\left(-\frac{\la^3}{6}+\frac{\mu a^2}{2}+\frac{a^3}{3}\right)
\cdot\sum_{r=0}^{\infty}f_r^La_1^{-3r}.
\endaligned\tag 2.3.26
$$
Now $a_1\sim a$, so by (2.2.23) and (2.2.18),
$$
\sum_{r=0}^{\infty}f_r^La_1^{-3r}\le_b\sum_{r\le L}\left(\frac{3r}{2e a_1^3}\right)^r+
\sum_{r>L}\left(\frac{3L}{2e a_1^3}\right)^r
\le\sum_{r\ge 0}\left(\frac{2}{e}\right)^r<\infty.
$$
So (2.3.26) becomes
$$
R_n^{(2)}\le_bn^{-1/12}\cdot n^{-1/3}e^{\ln^{3/4}n}
= n^{-5/12+o(1)}.
$$
Further, by (2.2.22) and $g_0^L=0$,
$$
\Sigma_2-\Sigma_1=N(n,p)\sum_{r>0}\left(\frac{p}{2q}\right)^r[x^n]\frac{g_r^Le^{H(x)}}
{(1-T(x))^{3r-1}}.
$$
Therefore
$$
\align
|R_n^{(1)}|\le_b&\, R_n^{(2)}+n^{-1/12}\cdot\bigl(n^{-1/3}\ln^{15/4}n\bigr)
\\
&\cdot \exp\left(-\frac{\la^3}{6}+\frac{\mu a^2}{2}+\frac{a^3}{3}\right)
\cdot\sum_{r=1}^{\infty}g_r^La_1^{-3r+1}.
\endalign
$$
So, using the bounds (2.2.19) and (2.2.23) for $g_r^L$, we conclude that 
$|R_n^{(1)}|\le 2R_n^{(2)}$. Thus, for $L=\ln^{1/4}n$,
$$
\Sigma_1^*+O(n^{-1/3+o(1)})\le \frac{\pr(S_n>0,\,\Cal E_n\le L)}
{(2\pi)^{1/2}e^{3/8}n^{-1/12}}\le \Sigma_2^*
+O(n^{-1/3+o(1)}).\tag 2.3.27
$$
where
$$
\aligned
\Sigma_2^*=&\,\sum_{r\ge 0}\left(\frac{np}{2q}\right)^rf_r^LA(1/4+3r,\mu),\\
\Sigma_1^*=&\Sigma_2^*- n^{-1/3}\sum_{r> 0}\left(\frac{np}{2q}\right)^rg_r^L
A(-3/4+3r,\mu).
\endaligned\tag 2.3.28
$$
\si

Let us have a close look at $\Sigma^*_1$ and $\Sigma_2^*$. Write
$$
\align
\Sigma_2^*=&\,\sum_{r\le L}\left(\frac{np}{2q}\right)^rf_rA(1/4+3r,\mu)+
\sum_{r>L}\left(\frac{np}{2q}\right)^rf_r^LA(1/4+3r,\mu)\\
=&\,\Sigma_{21}^*+\Sigma_{22}^*.
\endalign
$$
By $f_r^L\le f_r$, (2.2.18), (2.3.22), and Stirling's formula for $\Gamma(r)=(r-1)!$,
$$
\align
\Sigma_{22}^*\le_b&\,e^{2|\mu|^3/3}\sum_{r>L}\left(\frac{1}{2}+O\bigl(|\la|n^{-1/3}
\bigr)\right)^r
\left(\frac{3}{2}\right)^rr^{r-1/2}e^{-r}\left(\frac{2}{3}\right)^r\Gamma^{-1}(r)\\
\le&e^{2|\mu|^3/3}\sum_{r>L}\left(\frac{2}{3}\right)^r \le_b\left(\frac{2}{3}\right)^L .
\endalign
$$
Further, since uniformly  for $r\le L$,
$$
(1+\la n^{-1/3})^r=\exp\bigl(O(L|\la|n^{-1/3})\bigr)=1+O(n^{-1/3}\ln^{1/4} n),
$$
we have
$$
\Sigma_{21}^*=(1+o(1))\sum_{r\le L}\frac{f_r}{2^r}A(1/4+3r,\mu).
$$
And,  analogously to $\Sigma_{22}^*$,
$$
\sum_{r>L}\frac{f_r}{2^r}A(1/4+3r,\mu)\le b\left(\frac{2}{3}\right)^L.
$$
Therefore
$$
\Sigma_2^*\sim\sum_{r\le L}\frac{f_r}{2^r}A(1/4+3r,\mu)\to\sum_{r\le L}
\frac{f_r}{2^r}A(1/4+3r,\mu).
$$
Also, by the definition of $\Sigma_1^*$ in (2.3.25), it follows  that 
$|\Sigma_1^*-\Sigma_2^*|$
is $O(n^{-1/3})$. Hence
\proclaim{Proposition 2.3.2} For $|\la|=O(1)$,
$$
\frac{\pr(S_n>0,\,\Cal E_n\le L)}{(2\pi)^{1/2}e^{3/8}n^{-1/12}}
\sim c(\mu):=\sum_{r\ge 0}\frac{f_r}{2^r}A(1/4+3r,\mu),
$$
and $\mu\,(=\la+O(n^{-1/3}))$ can be replaced with $\la$, as $c(x)$ is
positive and continuous for all $x$.
\endproclaim
\bi

{\it Case $\la\to\infty$, $\la=o(n^{1/12})$.\/} According to (2.1.16),
we select $L=\a \la^3$, $\a>2/3$. This time we use a refined version of
(2.3.24), with the exponential factor sneaking behind the
sum operation for $r\le\a\la^3$, which allows us to choose $a\,(\le 2\la)$ dependent on $r$
for $r\le\a\la^3$.  Also, for those $r$ and $a$, $ra=O(\la^4)=o(n^{1/3})$ and $r\ln n=O(\la^3\ln n)
=o(n^{1/3})$; so $R_{n,3r}=\tilde\Delta_{n,3r}$ in this range.
For $r>\a\la^3$ we select $a=\la$, and here $R_{n,3r}=\Delta_{n,3r}$. (So, $a=o(n^{1/12})$
throughout.) By  (2.2.20),
$$
\align
&\exp\left(-\frac{\la^3}{6}+\frac{\la a^2}{2}+\frac{a^3}{3}\right)
\sum_{r>\a\la^3}(r+1)(1/2+o(1))^rf_r^L\\
\le_b&\exp\left(\frac{2\la^3}{3}\right)\sum_{r>\a\la^3}(r+1)\left(\frac{3\a\la^3(1+o(1))}
{4e\la^3}\right)^r\\
\le_b&\la^3\exp\left(\frac{2\la^3}{3}+\la^3\a\ln\frac{3\a(1+o(1))}{4e}\right)\\
\le_b&\la^3\exp\left(\la^3\a\ln\frac{3\a(1+o(1))}{4}\right),
\endalign
$$
and, pushing $\a$ down to $2/3$, we can make the coefficient of  $\la^3$
in the exponent arbitrarily close to
$$
\frac{2}{3}\cdot\ln\frac{1}{2}=-0.46\dots. 
$$
According to (2.2.18) and (2.3.13), it remains to bound
$$
 \sum_{r\le \a\la^3}\frac{r+1}{(r+1)^{1/2}}\min_{a\le 2\la}\left\{\exp\left(-\frac{\mu^3}{6}
+\frac{\mu a^2}{2}+\frac{a^3}{3}
\right)\left(\frac{3r}{4ea^3}\right)^r\right\};
$$
(we have replaced $a_1=n^{1/3}(1-e^{-an^{-1/3}})$ with $a$, since for $r\le \a\la^3$,
$$
a_1^{3r}=a^{3r}e^{O(\la^4 n^{-1/3})}\sim a^{3r},
$$
and $\la^3/6$ with $\mu^3/6$, see (2.3.18). So we need to find $\min_{a\le 2\la} \Phi(r,a)$,
$$
\Phi(r,a):=-\frac{\mu^3}{6}-r\ln\left(\frac{3r}{4ea^3}\right)+\frac{\mu a^2}{2}+
\frac{a^3}{3}.\tag 2.3.29
$$
$\Phi(r,a)$ attains its absolute minimum at $\xi(r)$, a unique positive root
of
$$
\mu \xi+\xi^2=\frac{3r}{\xi}\le \frac{3\a\la^3}{\xi},\tag 2.3.30
$$
i. e. $\xi(r)<2\la$ if $\a$ is sufficiently close to $2/3$ from above.
Further $\phi(r):=\Phi(r,\xi(r))$ attains its maximum at $\bar r$, a root
of
$$
\phi^\prime(y)=\ln\frac{3y}{4}-3\ln \xi(y)=0,
$$
i. e.
$$
\bar r=\frac{4}{81}\mu^3,\quad \bar a:=\xi(\bar r)=\frac{\mu}{3}\,(<2\la).\tag 2.3.31
$$
Consequently
$$
\phi(\bar r)=\exp\left(-\frac{10\mu^3}{81}\right)\sim\exp\left(-\frac{10\la^3}{81}\right).
\tag 2.3.32
$$
It is easy to show that
$$
\phi^{\prime\prime}(\bar r)=-\frac{27}{\mu^3}\sim -\frac{27}{\la^3},\tag 2.3.33
$$
and, with some work, that $\phi^{\prime\prime}(r)<0$ always. A standard application of Laplace method 
yields
$$
\sum_{r\le \a\la^3}(r+1)^{1/2}\min_{a\le\la}\left\{\exp\left(-\frac{\mu^3}{6}+
\frac{\la a^2}{2}+\frac{a^3}{3}
\right)\left(\frac{3r}{4ea^3}\right)^r\right\}
\le_b\la^3\exp\left(-\frac{10\la^3}{81}\right).
$$
Therefore, consulting (2.3.16) for $\tilde\Delta_{n,w}$ and (2.3.13) for $\Delta_{n,w}$,
we bound $R_n^{(2)}$, the total contribution of the remainders $(w+1)R_{n,w}$ to the sum
$\Sigma_{2}$ in (2.2.20):
$$
\align
R_n^{(2)}=\sum_r(r+1)R_{n,3r}\le_b&\,n^{-1/12}\exp\left(-\frac{10\la^3}{81}\right)
\bigl(n^{-1/3}\la^{3.75}
+n^{-1/3}\ln^{3/4}n\bigr)\\
&+n^{-1/12}\la^3e^{-0.27\la^3}
\bigl(n^{-1/3}\la^4+n^{-1/3}\ln^4 n\bigr)\\
\le_b&\,\bigl(\la^{-1/4}+n^{-1/3+o(1)}\bigr)n^{-1/12}
\exp\left(-\frac{10\la^3}{81}\right).\tag 2.3.34
\endalign
$$
As for  $R_n^{(1)}$, the total contribution of the remainders $(w+1)R_{n,w}$ to 
$\Sigma_1$ in (2.2.20),
it is $O(R_n^{(2)})$, just like the $\la=O(1)$ case.  So we arrive at the 
counterpart of (2.3.27)-(2.3.28),
with 
$$
\bigl(\la^{-1/4}+n^{-1/3+o(1)}\bigr)\exp\left(-\frac{10\la^3}{81}\right)
$$
 taking place of $n^{-1/3+o(1)}$.  Further, again we split
$\Sigma_2^*=\Sigma_{21}^*+\Sigma_{22}^*$. To bound $\Sigma_{22}^*$ we use the
second bound for $A(1/4+3r,\mu)$ with $a\equiv\lambda$ (2.3.21), and the bound (2.2.20) for
$f_r^L$. Just like $R_n^{(2)}$, we obtain
$$
\align
\Sigma_{22}^*\le_b&\,\la^{1/2}\exp\left(-\frac{\la^3}{6}+\frac{\la a^2}{2}+\frac{a^3}{3}
\right)\\
&\cdot\sum_{r>\a \la^3}(1/2+o(1))^rf_r^L\le_b\la^{1/2}e^{-0.46\la^3}.\tag 2.3.35
\endalign
$$
To evaluate sharply $\Sigma_{21}^*$, we use (2.3.23) to approximate
$A(y,\mu)$ for $\e \mu\le y$, $y=O(\la^3)$, and (2.3.24) to bound
$A(y,\mu)$ for $y\le\e\mu$, $\e>0$ sufficiently small. Invoking (2.2.18) as well,
we have
$$
\sum_{\e\mu\le r\le\a\la^3}\!\!\left(\frac{np}{2q}\right)^r\!\!f_rA(1/4+3r,\mu)
\sim\,
\frac{1}{4\pi}\!\sum_{\e\mu\le r\le\a\la^3}\!
\frac{\xi^{3/4} e^{\phi(r)}}{\bigl(r((3r+1)\xi^{-2}+\mu+2\xi)\bigr)^{1/2}};
$$
here $\phi(r):=\min_a \Phi(r,a)=\Phi(r,\xi)$, see (2.3.28)-(2.3.29) for $\Phi(r,a)$
and $\xi=\xi(r)$. We know that $\phi(r)$ attains its pronounced maximum at $\bar r=(4/27)\la^3$, 
i. e. well within $[\e\mu,\a\la^3]$. Using  (2.3.31)-(2.3.33), by Laplace method,
$$
\align
\sum_{\e\mu\le r\le\a\la^3}\!\!\left(\frac{np}{2q}\right)^r\!\!f_rA(1/4+3r,\mu)
\sim&\,\frac{1}{4\pi}\frac{\bar \xi^{3/4}}{\bigl(\bar r((3\bar r+1)\bar\xi^{-2}+\mu+2\bar\xi)\bigr)^{1/2}}
\left(\frac{2\pi}{-\phi^{\prime\prime}(\bar r)}\right)^{1/2}\\
\sim&\,\frac{1}{4(2\pi)^{1/2}3^{3/4}}\la^{1/4}\exp\left(-\frac{10\la^3}{81}\right).
\endalign
$$
Applying (2.3.24), it is not difficult to show that 
$$
\sum_{\ r\le\e\mu}\!\!\left(\frac{np}{2q}\right)^r\!\!f_rA(1/4+3r,\mu)\ll\la^{1/4}\exp\left(-\frac{10\la^3}{81}\right).
$$
So
$$
\Sigma_{12}^*:=\sum_{ r\le\a\la^3}\!\!\left(\frac{np}{2q}\right)^r\!\!f_rA(1/4+3r,\mu)
\sim\frac{1}{4(2\pi)^{1/2}3^{3/4}}\la^{1/4}\exp\left(-\frac{10\la^3}{81}\right),
$$
hence (see (2.3.35))
$$
\Sigma_2^*\sim\frac{1}{4(2\pi)^{1/2}3^{3/4}}\la^{1/4}
\exp\left(-\frac{10\la^3}{81}\right),\tag 2.3.36
$$
as well. And, analogously to the $\la=O(1)$ case, for $\Sigma_1^*$ defined in (2.3.28),
$$
\bigl|\Sigma_1^*-\Sigma_2^*|\ll \la^{1/4}\exp\left(-\frac{10\la^3}{81}\right).
\tag 2.3.37
$$
\proclaim{Proposition 2.3.3} For $\la\to\infty$, $\la=o(n^{1/12})$,
$$
\pr(S_n>0,\,\Cal E_n\le L)\sim\frac{e^{3/8}}{4\cdot 3^{3/4}}\,\,\la^{1/4}\exp\left(-\frac{10\la^3}{81}
\right).
$$
\endproclaim
\si

{\bf Proof of Proposition 2.3.3\/} The probability is asymptotic to the 
expression in (2.3.36)
times $(2\pi)^{1/2}e^{3/8}n^{-1/12}$.\qed
\bi

Lastly,
\si
{\it Case $\la\to-\infty$, $|\la|=o(n^{1/12})$.\/}  According to (2.1.16), we can pick  $L=0$. 
By Proposition 2.2.4 and  (2.3.17) for $w=0$, and $a\ge |\la|$, $a=o(n^{1/12})$, we have
$$
\align
\pr(S_n>0,\Cal E_n\le 0)=&\,N(n,p)\,[x^n]\,e^{H(x)}\\
=&\,(2\pi)^{1/2}n^{-1/12}e^{3/8}A(1/4,\mu)+O(\Delta_{n,0}).
\endalign
$$
Notice that 
$$
\left.\left(\frac{\mu a^2}{2}+\frac{\a^3}{3}\right)\right|_{a=|\la|}=\frac{\la^3}{6}+o(1),
$$
since $\la^3-\mu^3=o(1)$.  Setting $a=\la$ in (2.3.13), we obtain
$$
\Delta_{n,0}\ll  n^{-1/12}\cdot n^{-1/3}\cdot n^{-1/3}\bigl(|\la|^{3.75}+\ln^{3.75} n\bigr).
$$
And, by (2.3.21),
$$
A(1/4,\mu)\sim(2\pi)^{-1/2}|\mu|^{1/2-1/4}\sim (2\pi)^{-1/2}|\la|^{1/4}.
$$
\proclaim{Proposition 2.3.4} Suppose $\la\to -\infty$, $|\la|=o(n^{1/12})$. Then
$$
\pr(S_n>0,\Cal E_n\le 0)\sim e^{3/8}|\la|^{1/4}.
$$
\endproclaim
\bi

Since in each of the three cases our $L$ is such that
$$
\lim\pr(\Cal E_n\le L)=1,
$$
Propositions 2.3.2 -2.3.4  combined with the relations (2.1.12) and (2.1.13), prove
the part of Theorem 1.1 about $G(n,p)$, $\hat p=1/2$.
\bi
{\bf 3. Solvability probability: $G(n,m)$ and $\hat p=1/2$.\/}
\bi

Our task is to show that the result for the near-critical $G(n,p)$, $p=(1+\la n^{-1/3})/n$,
$\la=o(n^{1/12})$, implies the analogous claim for $G(n,m)$, $m=(n/2)(1+\la n^{-1/3})$.
Denoting $N=\binom{n}{2}$, 
$$
p=\frac{m}{N}+O(m^{1/2}N^{-1})=\frac{1+n^{-1/3}\la^\prime}{n},\quad\la^\prime=\la+O(n^{-1/6})).
\tag 3.1
$$
Obviously $\la^\prime=o(n^{1/12})$, so 
$$
\pr(S(G(n,p)>0)\to 0.
$$
Since an event $\{S(G)>0\}$ is monotone (increasing) with $G$, a general ``$p$-to-$m$''
result, Bollob\'as [5], \L uczak [15], implies that
$$
\pr(S(G(n,m)>0)\to 0,
$$
too. However we want to prove a sharp formula
$$
\pr(S(G(n,m)>0)\sim c(\la)n^{-1/12},\tag 3.2
$$
so that the probabilities in question can be as small as
$$
\exp\left(-\frac{10}{81}(n^{1/12-o(1)})^3\right)=\exp\left(-\frac{10}{81}n^{1/4-o(1)}\right).
$$
It turns out that in our case the argument in [5], [15] can be sharpened to yield (3.2).
\si

To start, recall the classic entropy bound
$$
\align
\pr(\text{Bin}(N,p)\ge k)\le&\, \exp[NH(k/N)],\quad k> Np,\\
\pr(\text{Bin}(N,p)\le k)\le&\, \exp[NH(k/N)],\quad k< Np,
\endalign
$$
where 
$$
H(x):x\ln (p/x)+(1-x)\ln(q/(1-x)).
$$
Approximating $H(x)$ by its second degree Taylor polynomial plus a remainder term, we obtain:
uniformly for  $p\le 1/2$, and $\omega\le a(Np)^{1/6}$, $a>0$ being fixed,
$$
\aligned
\pr\left(\text{Bin}(N,p)\ge Np+\omega\sqrt{Npq}\right)\le_b&\,e^{-\omega^2/2},\\
\pr\left(\text{Bin}(N,p)\le Np-\omega\sqrt{Npq}\right)\le_b&\,e^{-\omega^2/2},
\endaligned\tag 3.3
$$
(The bounded factor implicit in $\le_b$ notation depends on $a$.) Given $m$ and $\omega\le
m^{1/6}$, introduce  $p_1<p_2$:
$$
\aligned
Np_1+\omega \sqrt{Np_1}=&m\,\,\Longrightarrow p_1=(4N)^{-1}\bigl(\sqrt{4m+\omega^2}-\omega)^2,\\
Np_2-\omega \sqrt{Np_2}=&m\,\,\Longrightarrow p_2=(4N)^{-1}\bigl(\sqrt{4m+\omega^2}+\omega)^2.
\endaligned\tag 3.4
$$
Then
$$
\aligned
\frac{Np_2}{\omega^6}>\frac{Np_1}{\omega^6}=&\,\frac{m}{\omega^6}\left(\sqrt{1+\omega^2/4m}-
\omega/(2\sqrt{m}\right)^2\\
\ge&\,a:=(\sqrt{2}-1)^2,
\endaligned\tag 3.5
$$
as $\omega/2\sqrt{m}\le 0.5m^{-1/3}\le 1$.
\si

Now, using $e(G)$ to denote the number of edges in a graph $G$,  $e(G(n,p))=\text{Bin}(N,p)$. So,
by  (3.3)-(3.5), 
$$
\aligned
&\pr(e(G(n,p_1))>m)\le \pr(e(G(n,p_1))\ge Np_1+\omega\sqrt{Np_1q_1})\le_be^{-\omega^2/2},\\
&\pr(e(G(n,p_2))< m)\le \pr(e(G(n,p_2))\le Np_2-\omega\sqrt{Np_2q_2})\le_be^{-\omega^2/2}.
\endaligned\tag 3.6
$$
Since
$$
\pr(S(G(n,p))>0)=
\sum_{\mu=0}^N\pr(e(G(n,p))=\mu)\pr(S(G(n,\mu))>0),
$$
and $\pr(S(G(n,\mu))>0)$ decreases with $\mu$, we have
$$
\align
\pr(S(G(n,p_1))>0)\ge& \pr(e(G(n,p_1))\le m)\pr(S(G(n,m))>0)\\
\ge& (1-O(e^{-\omega^2/2}))\pr(S(G(n,m))>0),
\endalign
$$
and
$$
\align
\pr(S(G(n,p_2))=0)\ge&\pr(e(G(n,p_2))\ge m)\pr(S(G(n,m))=0))\\
=& (1-O(e^{-\omega^2/2}))\pr(S(G(n,m))=0).
\endalign
$$
Therefore
$$
\multline
 \frac{\pr(S(G(n,p_1))>0)}{1-O(e^{-\omega^2})}\ge\pr(S(G(n,m))>0)\\
 \ge
 \frac{\pr(S(G(n,p_2))>0)-O(e^{-\omega^2/2})}{1-O(e^{-\omega^2/2})}.
 \endmultline\tag 3.7
$$
Now, by (3.4),
$$
\align
p_{1,2}=&\frac{m}{N}\bigl(1+O(\omega m^{-1/2})\bigr)
=\frac{1+\la_{1,2}n^{-1/3}}{n},\\
\la_{1,2}=&\la +O(\omega m^{-1/2}+n^{-2/3}),
\endalign
$$
so, as $|\la|=o(n^{1/12})$, 
$$
\la_{1,2}^3=\la^3+O\bigl[\la^2(\omega m^{-1/2}+n^{-2/3})\bigr]+
O\bigl[(\omega m^{-1/2})^3+n^{-2}\bigr]=\la^3+o(\omega n^{-1/3}).
$$
That is, $\la_{1,2}^3-\la^3\to 0$. Hence, 
$$
\pr(S(G(n,p_{1,2}))>0)\sim c(\la)n^{-1/12}.\tag 3.8
$$
Also $\omega^2\gg |\la|^3$ if $\omega=n^{1/8}$, which is compatible with
the restriction $\omega\le n^{1/6}$. For this choice of $\omega$, the
relations (3.7)-(3.8) imply: for $\la=o(n^{1/12})$,
$$
\pr(S(G(n,m))>0)\sim c(\la)n^{-1/12}.
$$
This completes the proof of Theorem 1.1 for $\hat p=1/2$.\qed
\bi
{\bf 4. Solvability ($2$-colorability) probability: $G(n,p)$, $G(n,m)$ and $\hat p=1$.\/}
\si

Consider the $G(n,p)$ case. We know that the system
$$
x_i+x_j\equiv 1\,(\text{mod }2),\quad (i,j)\in E(G)
$$
is solvable iff the graph $G$ has no odd cycles. So a counterpart of 
(2.1.11) is
$$
\pr(S_n>0,\,\Cal E_n\le L)=\,N(n,p)\,[x^n]\,
\exp\left[\sum_{\ell= -1}^L\left(\frac{p}{q}\right)^{\ell}C_{\ell}^e(x)
\right],\tag 4.1
$$
($S_n=S(G(n,p))$, $\Cal E_n=\Cal E(G(n,p))$), where $C_{\ell}^e(x)$ is the 
exponential generating function of graphs $G$ {\it without odd
cycles\/}, with an excess $\Cal E(G)=\ell$. And again the events
$\{S_n>0\}$ and $\{\Cal E_n\le L\}$ are positively
correlated, i. e.
$$
\pr(S_n>0,\,\Cal E_n\le L)\le \pr(S_n>0)\le
\frac{\pr(S_n>0,\,\Cal E_n\le L)}{\pr(\Cal E_n\le L)}.
$$
Thus the generating functions $C_{\ell}^e(x)$ take a center stage. Obviously
$$
C_{-1}^e(x)=C_{-1}(x)\left(=T(x)-\frac{1}{2}\,T^2(x)\right).
$$
Furthermore, while
$$
C_0(x)=\frac{1}{2}\left(\ln\frac{1}{1-T(x)}-T(x)-\frac{1}{2}\,T^2(x)\right),
$$
for $C_0^e(x)$ we have
$$
C_0^e(x)=\frac{1}{4}\left(\ln\frac{1}{1-T^2(x)}-T^2(x)\right).\tag 4.2
$$
Indeed, we enumerate the connected 
unicyclic graphs with an even cycle, i. e. forests of an even number of rooted 
trees,
whose roots form an undirected cycle. So
$$
C_0^e(x)=\sum_{\text{ even }j\ge 4}\frac{(j-1)!}{2}\frac{T^j(x)}{j!},
$$
which simplifies to (4.2). Comparing $C_0^e(c)$ and $C_0(x)$ we see that,
for $|x|<e^{-1}$, $x\to 1$, i. e. for $x$ dominant asymptotically,
$$
C_0^e(x)=\frac{1}{2}C_0(x)+\frac{1}{8}-\frac{1}{4}\ln 2+O(|T(x)-1|);\tag 4.3
$$
in particular, $C_0^e(x)\sim (1/2)C_0(x)$. 
We want to show that this pattern persists for $\ell>0$, namely 
$$
C_{\ell}^e(x)\sim \frac{1}{2^{\ell+1}}C_{\ell}(x),\quad (|x|<e^{-1},\,x\to e^{-1}).
\tag 4.4
$$
Comparing  (2.1.11) and  (4.1), and recalling the different roles played
by $C_0(x)$ and $\{C_{\ell}(x)\}_{\ell>0}$ in the analysis of the $\hat p=1/2$
makes it transparent, hopefully, that for $\hat p=1$  
we should have
$$
\pr(G(n,p)\text{ is }2\text{-colorable})=\pr(S_n>0)\sim 2^{-1/4}e^{1/8}
c(\la)n^{-1/12}.
$$
\si

Let us prove (4.4). First
\proclaim{Proposition 4.1} Given $n$ and $m\le N:=\binom{n}{2}$, let $C(n,m)$ denote the total
number of connected graphs on $[n]$ with $m$ edges, and let $C^e(n,m)$ denote the
total number of connected graphs without odd cycles. Then
$$
C^e(n,m)\le \frac{1}{2^{m+1-n}}\,C(n,m).\tag 4.5
$$
Consequently
$$
C_{\ell}^e(x)\le_c\frac{1}{2^{\ell+1}}C_{\ell}(x),\quad\ell\ge -1.\tag 4.6
$$
\endproclaim
\si

{\bf Proof of Proposition 4.1.\/} We begin with a simple claim.
\si
\proclaim{Lemma 4.2} Let $T$ be a tree on the vertex set $[n]$. Let $X(T)$ denote
the total number of paths in $T$ of an even edge-length $2$ at least. Then 
$X(T)\ge X(P_n)$, where $P_n$ is a path on $[n]$, and 
$$
X(P_n)=\left\lceil \frac{n(n-2)}{4}\right\rceil. \tag 4.7
$$
\endproclaim
\si

{\bf Proof of Lemma 4.2.\/} Pick a vertex $v\in [n]$, and introduce $V_0(T)$ and $V_1(T)$
the set of vertices reachable from $v$ by paths of even length $2$ at least, and 
odd length
respectively; in particular $v\in V_0$. Now every two vertices from $V_i(T)$, 
$(i=0,1)$,
are connected by an even path, while there is no even path connecting 
$v_0\in V_0(T)$ 
and $v_1\in V_1(T)$.
Hence
$$
X(T)=\binom {|V_0(T)|}{2}+\binom {|V_1(T)|}{2}.
$$
It follows that $X(T)$ attains its minimum when $|V_0(T)|=\lfloor n/2\rfloor$ 
and
$|V_1(T)|=\lceil n/2\rceil$, or the other way around, i. e. when $T=P_n$, and the
minimum value is
$$
X(P_n)=\binom{\lfloor n/2\rfloor}{2}+\binom{\lceil n/2\rceil}{2}
=\left\lceil \frac{n(n-2)}{4}\right\rceil.
$$
\qed
\si

Armed with this Lemma, we will derive a recurrence {\it inequality\/}
for $C^e(n,m)$. First we recall a recurrence equality for
$C(n,m)$, [23], [3]: for $n\ge 3$, $n-1\le m\le N$,
$$
\multline
mC(n,m)=(N-m+1)C(n,m-1)\\
+\frac{1}{2}\sum\limits_{n_1+n_2=n,\atop m_1+m_2=m-1}\binom{n}{n_1}
n_1n_2C(n_1,m_1)C(n_2,m_2).
\endmultline\tag 4.8
$$
Explanation. The left hand side of (4.8) is the total number
of the connected $(n,m)$ graphs with a marked edge. Each one of these graphs with 
a marked edge can be obtained in one of two, mutually exclusive ways. First way is
inserting a marked edge into a connected graph on $[n]$ with $m-1$ edges, which accounts 
for the first term on the right hand side of (4.6); indeed $N-m+1$ is the total number
of unordered pairs of vertices not connected by an edge in a given connected graph
with $m-1$ edges. Second way is to start with a connected $(n_1,m_1)$ graph and
a connected $(n_2,m_2)$ graph, having $m-1$ edges in total, and to add a marked edge
that joins two connected graphs; $n_1n_2$ is the total number of ways to select two
``contact'' points, representing each of two graphs. 
\si

Let us see if there is a similar recursive formula for $C^e(n,m)$. Clearly, if a
marked edge joins two connected graphs, none of these two graphs may have an odd
cycle. So we definitely have the ``$C^e(\cdot,\cdot)$'' counterpart of the second
term on the right hand side of (4.8). As for a potential counterpart of the first term,
a difficulty is that an additional $m$-th edge is not allowed to form an
odd cycle with any of the $m-1$ edges already present. And so the total number of
admissible options depends on the structure of a $(n,m-1)$ graph $G$ in question.
(For such a graph to be connected, it is necessary that $m\ge n$.) However we can bound 
the number of options.  $G$ is spanned by a 
tree $T$ on $[n]$, and none of the $m-1-(n-1)=m-n$ edges of $G\setminus T$ completes
an odd cycle by joining the ends of an even path in $T$. By Lemma 4.2, the total number of 
those even paths is $\lceil n(n-2)/4\rceil$, at least. Hence the total number
of options for the $m$-edge is $N-(m-1)-\lceil n(n-2)/4\rceil$, at most. And it is 
straightforward that, for $n\ge 3$ and by $m\ge n$,
$$
N-(m-1)-\left\lceil\frac{n(n-2)}{4}\right\rceil\le \frac{1}{2}(N-(m-1)).
$$ 
So 
$C^e(\cdot,\cdot)$ satisfies a recursive {\it inequality\/}: for $n\ge 3$, $n-1\le m\le N$,
$$
\multline
mC^e(n,m)\le\frac{1}{2}(N-m+1)C^e(n,m-1)\\
+\frac{1}{2}\sum\limits_{n_1+n_2=n,\atop m_1+m_2=m-1}\binom{n}{n_1}
n_1n_2C^e(n_1,m_1)C^e(n_2,m_2).
\endmultline\tag 4.9
$$
($C^e(\nu,\mu):=0$ if $\nu=0$, or $\mu\notin [\nu-1,\binom{\nu}{2}]$.) We will use
(4.9) and induction to prove the bound (4.5).
To this end, we define a lexicographical order, $\prec$,
on $\{(n,m)\,:\,n\ge 1,\,n-1\le m\le \binom{n}{2}\}$ as follows: denoting $\ell=m-n$,
$$
(n_1,m_1)\prec (n_2,m_2) \Longleftrightarrow \ell_1<\ell_2,\text{ or }\ell_1=\ell_2
\text{ and }n_1<n_2.
$$
The order $\prec$ is total, and $(1,0)$ is the minimal element. 
The inductive basis holds, since $C^e(1,0)=C(1,0)=1$, and  $C^e(2,1)=C(2,1)=1$.  Suppose 
that, for some $n\ge 2$ and $m\in [n-1,\binom{n}{2}]$,
$$
C^e(\nu,\mu)\le \frac{1}{2^{\mu-\nu+1}}C(\nu,\mu),\quad \forall\, (\nu,\mu)\prec (n,m).
$$
Since $(n,m-1)\prec (n,m)$, the inductive assumption implies that
$$
\align
\frac{1}{2}(N-m+1)C^e(n,m-1)\le&\, \frac{1}{2}(N-m+1)\frac{1}{2^{m-1-n+1}}C(n,m-1)\\
=\,&\frac{1}{2^{m-n+1}}(N-m+1)C(n,m).\tag 4.10
\endalign
$$
Further, for the double sum in (4.9),
$$
m_1-n_1+1\ge 0,\quad m_2-n_2+1\ge 0,
$$
and
$$
(m_1-n_1+1)+(m_2-n_2+1)=m-1-n+2=m-n+1,
$$
so that
$$
m_i-n_i+1\le m-n +1\Longrightarrow m_i-n_i\le m-n,\quad i=1,2.
$$
So, for $n_1,\,n_2>0$, we have $(n_i,m_i)\prec (n,m)$ and therefore, by the
inductive assumption,
$$
\prod_{i=1}^2C^e(n_i,m_i)\le\prod_{i=1}^2\frac{1}{2^{m_i-n_i+1}}C(n_i,m_i)
=\frac{1}{2^{m-n+1}}\prod_{i=1}^2C(n_i,m_i).\tag 4.11
$$
Combining (4.9)-(4.11), and the recurrence equation (4.8) for $C(\cdot,\cdot)$, we 
obtain
$$
\align
mC^e(n,m)\le&\,\frac{1}{2^{m-n+1}}(N-m+1)C(n,m-1)\\
&+\frac{1}{2^{m-n+1}}\,\frac{1}{2}\sum\limits_{n_1+n_2=n,\atop m_1+m_2=m-1}\binom{n}{n_1}
\prod_{i=1}^2n_iC(n_i,m_i)\\
=&\,\frac{1}{2^{m-n+1}}\,mC(n,m).
\endalign
$$
Thus the bound (4.3) holds for $(n,m)$ too. The proof of Proposition
4.1 is complete.\qed
\si

By Proposition 4.1 and and the formula (4.1), we have
$$
\multline
\pr(S_n>0,\,\Cal E_n\le L)\\
\le \,N(n,p)\,[x^n]\,
\exp\left[\frac{q}{p}C_1(x)+C_0^e(x)+\frac{1}{2}\sum_{\ell= 1}^L\left(\frac{p}{2q}
\right)^{\ell}C_{\ell}(x)\right].\endmultline\tag 4.12
$$
Since $C_0^e(x)$ is asymptotic to $(1/2)C_0(x)+\ln(2^{-1/4}e^{1/8})$ as
$x\to e^{-1}$, only a trivial change in the proof of
Theorem 1.1 {\bf (i)\/} is needed to show that
$$
\pr(S_n>0)\lesssim 2^{-1/4}e^{1/8}c(\la)n^{-1/12},\quad (|\la|=o(n^{1/12}).\tag 4.13
$$
We omit the details.
Furthermore, since for $\la\to -\infty$ we use $L=0$, the sums $\sum_{\ell=1}^L$ in
(4.1), (4.12) disappear, and we obtain an asymptotic equality
$$
\pr(S_n>0)\sim 2^{-1/4}e^{1/8}c(\la)n^{-1/12},\quad (|\la|=o(n^{1/12},\,\la\to 
-\infty).
$$
To complete the proof of {\bf (ii)\/}, (case $\la=O(1)$), we need to prove (4.4)  
for each {\it fixed\/} $\ell>0$. Recall Wright's formula
$$
C_{\ell}(x)=(1-T(x))^{-3\ell}\left[\sum_{d=0}^{2\ell}c_{\ell,d}(1-T(x))^d\right],
\quad (\ell>0),\tag 4.14
$$
Let us find  a similar formula for $C_{\ell}^e(x)$, $\ell> 0$.
\proclaim{Proposition 4.3}  For $\ell>0$,
$$
C^e_{\ell}(x)=(1-T^2(x))^{-3\ell}\left[\sum_{d=0}^{8\ell-1}c^e_{\ell,d}
(1-T(x))^d\right],
\tag 4.15
$$
where
$$
c^e_{\ell,0}=2^{2\ell-1}\,c_{\ell,0}.\tag 4.16
$$
Consequently, for $|x|<e^{-1}$ and $x\to e^{-1}$, 
$$
C^e_{\ell}(x)\sim\frac{1}{2^{\ell+1}}C_{\ell}(x).
$$
\endproclaim 
{\bf Proof of Proposition 4.3.\/} We use the ideas of Wright's original proof of
(4.14), and the improvements suggested by Stepanov [22], (cf.  
[12], Section 9).
\si

Given a connected graph $G$ on $[n]$, with an excess $\ell=e(G)-v(G)>0$, we apply a
``pruning'' algorithm which successively deletes vertices of degree $1$. Obviously
the excess is preserved, and so for a terminal graph (core) $\bar G$ we have
$e(\bar G)-v(\bar G)=\ell$. $\bar G$ inherits all the cycles of $G$, and thus $\bar G$
has only even cycles iff $G$ does. A minimum degree of $\bar G$ is
$2$ at least, and---since $\ell(\bar G)=\ell>0$---a maximum degree is $3$ at least. Next we apply a ``cancellation'' 
algorithm to $\bar G$: at each step, we delete a vertex of degree $2$, splicing 
together the two edges it formerly touched. The excess is preserved again. Once all 
the vertices of degree $2$ are
gone, we get a connected multigraph (kernel) $\tilde G$, with possible loops and 
parallel edges, and a minimum vertex degree $3$ at least. Thus
$$
2e(\tilde G)\ge 3v(\tilde G),\quad e(\tilde G)-v(\tilde G)=\ell,
$$
and so
$$
v(\tilde G)\le 2\ell, \quad e(\tilde G)\le 3\ell. \tag 4.17
$$
Notice that the largest numbers of vertices and the edges in the kernel are
$2\ell$ and $3\ell$ respectively, and the corresponding kernel is a
$3$-regular multigraph. (In [12] graphs $G$ with such kernels were called
clean. It is these clean graphs that are most populous asymptotically
among all connected graphs on $[n]$ with excess $\ell$.)
Now that we have a reduced number $v(\tilde G)$ of vertices, we relabel them
using indices from $[v(\tilde G)]$ and preserving the order of their old indices from
$[n]$. Under this rule, it follows from (11) that the number of kernels $\tilde G$
for the collection of all connected graphs $G$ on $[n]$ with excess $\ell$ is a function 
of $\ell$ only!
\si

A key element of Wright's argument was the following identity. Let $M$ be a
connected multigraph on a vertex set $[\nu]$, with $\mu_i$ indistinguishable
loops at vertex $i$, and $\mu_{ij}$ indistinguishable parallel edges joining $i$ and
$j$, ($i,j\in [\nu]$, $i\neq j$). Let $h_{n,M}$ denote the total number
of the connected simple graphs $G$ on $[n]$, with minimum degree $2$ at
least and maximum degree $3$ at least (core-type graphs, in short), such that 
$\tilde G=M$. Letting
$$
\Cal H_M(z)=\sum_n\frac{h_{n,M}}{n!}\,z^n,\tag 4.18
$$
we have
$$
\Cal H_M(z)=\frac{\kappa}{\nu!}\,\frac{z^{\nu}}{(1-z)^{\mu}}\cdot\Cal K_M(z),\quad
\mu:=\sum_i\mu_i+\sum_{i<j}\mu_{ij},\tag 4.19
$$
where
$$
\aligned
\Cal K_M(z)=&\prod_{1\le i\le\nu}\left(z^{2\mu_i}\prod_{1\le i<j\le \nu}
z^{\mu_{ij}-1}(\mu_{ij}-(\mu_{ij}-1)z)\right),\\
\kappa=&\prod_{1\le i\le \nu}\frac{1}{2^{\mu_i}\mu_i!}\prod_{1\le i<j\le \nu}
\frac{1}{\mu_{ij}!}.
\endaligned\tag 4.20
$$
(Observe that $\Cal K_M(1)=1$.) Once (4.19) is established, it is easy to determine 
$H_M(x)$, 
the exponential generating 
function of all connected graphs $G$ whose kernel is the multigraph $M$. Indeed to go
from a core $\bar G$ back to $G$ on $[n]$  we need to choose an ordered sequence of
$v(\bar G)$ of rooted trees, of total size $n$, and plant them at the vertices of 
$\bar G$, moving increasingly from vertex $1$ to vertex $v(\bar G)$. Since
the generating function of such sequences is $T(x)^{v(\bar G)}$, we see that
$$
H_M(x)=\sum_n\frac{h_{n,M}}{n!}\,T(x)^n=\Cal H_M(T(x)).\tag 4.21
$$
Finally
$$
\aligned
C_{\ell}(x)=&\sum_{M:\,e(M)-v(M)=\ell}\!\!\!\!H_M(x)
=(1-T(x))^{-3\ell}\left[\sum_{d=0}^{2\ell}c_{\ell,d}(1-T(x))^d\right];\\
c_{\ell,0}:=&\frac{1}{(2\ell)!}\sum_{\boldsymbol\mu\text{ meets }(4.18)}\!\!\!
\kappa(\boldsymbol\mu).
\endaligned\tag 4.22
$$
\si

Our first step is to obtain a counterpart of (4.18)-(4.20) for 
$$
\Cal H_M^e(z)=\sum_n\frac{h^e_{n,M}}{n!}\,z^n,
$$
where $h^e_{n,M}$ is the total number of the connected core-type graphs on $[n]$
with only even cycles, that cancel to a given multigraph $M$. To this end, consider
an auxilliary problem. Let
$$
\mu_i=\mu_i^e+\mu_i^o,\quad \mu_{ij}=\mu_{ij}^e+\mu_{ij}^o,\quad (1\le i\neq j\le \nu).
\tag 4.23
$$
Let $h_{n,(\boldsymbol{\mu}^e,\boldsymbol{\mu}^o)}$ denote the total number of the
core-type graphs $G$ on $[n]$, which cancel to $M$, such that: (1) for each $i$, $G$ has an
even (odd resp.) number of $2$-degree vertices put on each of $\mu_i^e$ 
($\mu_i^0$ resp.) loops at vertex $i$ of $M$; (2) for each $(i,j)$, $G$ has an even (odd
resp.) number of $2$-degree vertices put on each of $\mu_{ij}^e$ ($\mu_{ij}^o$ resp.)
parallel edges joining the vertices $i$ and $j$ in $M$. Let us determine
$$
\Cal H_{(\boldsymbol{\mu}^e,\boldsymbol{\mu}^o)}(z)=\sum_n\frac{h_{n,
(\boldsymbol{\mu}^e,\boldsymbol{\mu}^o)}}{n!}\,z^n.
$$
A core-type graph $G$ on $[n]$ cancelling to $M$ and meeting the parity
conditions (1)-(2) can be viewed as a partition of $[n]$ into:
\si
(a) a subset of cardinality $\nu$, whose elements are the assigned to the $\nu$
vertices of $M$ in a unique (order-preserving) fashion;
\si
(b) $\forall\,i\in [\nu]$, a collection of $\mu_i$ {\it ordered\/} subsets, each
having $2$ elements at least (as $G$ is simple), such that exactly $\mu_i^e$
($\mu_i^o$ resp.) subsets have an even (odd resp.) number of elements;
\si
(c) $\forall\, 1\le i\neq j\le \nu$, a collection of $\mu_{ij}$ ordered
subsets, with at most one empty subset (as $G$ is simple), such that exactly 
$\mu_{ij}^e$ ($\mu_{ij}^o$ resp.) subsets have an even (odd resp.) number of
elements.
\si

So $\Cal H_{(\boldsymbol{\mu}^e,\boldsymbol{\mu}^o)}(z)$ is the product of 
generating functions $\Cal H_t(z)$ corresponding to $1+\nu+\binom{\nu}{2}$
combinatorial structures described in (a), (b), (c). The first is easy:
$$
\Cal H(z)=\frac{z^{\nu}}{\nu!}.\tag 4.24
$$
Next, for $i\in [\nu]$,
$$
\Cal H_{i}(z)=\frac{1}{2^{\mu_i}\mu_i^e!\mu_i^o!}\sum_n a_{n, (\mu_i^e,\mu_i^o)}z^n;
$$
here $a_{n, (\mu_i^e,\mu_i^o)}$ is the total number of {\it compositions\/} of
$n$ with $\mu_i$ parts, each $2$ at least, such that the first $\mu_i^e$ parts
(the last $\mu_i^0$ parts resp.) are even (odd resp.). The factor $1/2^{\mu_i}$ is 
needed as we do not distinguish between two opposite orderings of vertices sprinkled
on each of $\mu_i$ loops of $M$ at $i$. Consequently
$$
\aligned
\Cal H_{i}(z)=&\frac{1}{2^{\mu_i}\mu_i^e!\mu_i^o!}
\left(\sum_{k\ge 1}z^{2k}\right)^{\mu_i^e}\left(\sum_{k\ge 1}z^{2k+1}\right)^{\mu_i^o}
\\
=&\frac{1}{2^{\mu_i}\mu_i^e!\mu_i^o!}\frac{z^{2\mu_i+\mu_i^o}}{(1-z^2)^{\mu_i}}.
\endaligned\tag 4.25
$$
Similarly, for each $1\le i<j\le \nu$,
$$
\Cal H_{ij}(z)=\frac{1}{\mu_{ij}^e!\mu_{ij}^o!}\left[\left(\frac{z^2}{1-z^2}
\right)^{\mu_{ij}^e}+\mu_{ij}^e\left(\frac{z^2}{1-z^2}
\right)^{\mu_{ij}^e-1}\right]\left(\frac{z}{1-z^2}\right)^{\mu_{ij}^o}.
\tag 4.26
$$
Taking the product of the generating functions in (4.24)-(4.26) we obtain
$$
\Cal H_{(\boldsymbol{\mu}^e,\boldsymbol{\mu}^o)}(z)=
\frac{\kappa}{\nu!}\frac{z^{\nu}}{(1-z^2)^{\mu}}\cdot 
\Cal K_{(\boldsymbol{\mu}^e,\boldsymbol{\mu}^o)}(z), \tag 4.27
$$
where
$$
\multline
\Cal K_{(\boldsymbol{\mu}^e,\boldsymbol{\mu}^o)}(z)\\
=\prod_{i\in [\nu]}\binom{\mu_i}{\mu_i^e}z^{2\mu_i+\mu_i^o}
\cdot\prod_{1\le i<j\le \nu}\binom{\mu_{ij}}{\mu_{ij}^e}z^{\mu_{ij}+\mu_{ij}^e-2}
\left[\mu_{ij}^e-(\mu_{ij}^e-1)z^2\right].
\endmultline\tag 4.28
$$
As a partial check, summing over $(\boldsymbol{\mu}^e,\boldsymbol{\mu}^o)$, we
obtain Wright's formula (4.18)-(4.20).
\si

Now, for a core-type graph $G$ on $[n]$ without odd cycles, 
that cancels to $M$, $G$'s parity parameters $\boldsymbol{\mu}^e$, $\boldsymbol{\mu}^o$
must satisfy certain conditions.
First of all, for each $i\in[\nu]$, $\mu_i^e=0$,
since otherwise $G$ would have an odd cycle, with a single branching vertex. 
Likewise, for $\mu_{ij}>0$, the numbers of 
$2$-degree vertices of $G$ on $\mu_{ij}$ parallel edges of $M$ must all be of
the same parity, hence $\mu_{ij}^e=\mu_{ij}$ or $\mu_{ij}^o=\mu_{ij}$.
Subject to this condition, how many choices for 
$(\boldsymbol{\mu}^e,\boldsymbol{\mu}^o)$ do we have? For each $(i,j)$ such that
$\mu_{ij}>0$, define
$$
b_{ij}=b_{ji}=\left\{\alignedat2
&1,\quad&&\text{if }\mu_{ij}^e=\mu_{ij},\\
&0,\quad&&\text{if }\mu_{ij}^e=0.\endalignedat\right.
$$
If $C$ is a cycle in $M$, then the parity of a cycle in $G$ that cancels to $C$ is
the parity of $b(C):=\sum_{(i,j)\in C}b_{ij}$. Hence $b(C)$ must be even for
all cycles $C$, and we need to check this condition only for simple cycles that
do not use parallel edges. Let $T=T(M)$ be a tree on $[\nu]$ that spans $M$. 
Pick $\mu_{ij}^e$ for all $\nu-1$ pairs $(i,j)$ such that $(i,j)\in E(T)$, i. e.
one of $\mu_{ij}$ parallel edges is in $E(T)$. Let $\mu_{ij}>0$ and $e=(i,j)\notin
E(T)$. Then $e$ completes a cycle $C$ with a path in $T$ that connects $i$ and
$j$. The condition ``$b(C)$ is even'' determines $\mu_{ij}^e$ uniquely. Hence
a choice of $\nu-1$ values of $\mu_{ij}^e$ determines uniquely the remaining
$\mu^e_{\cdot\cdot}$. Arguing as in the proof of Lemma 2.1.1, we see that the
condition ``$b(C)$ is even'' will hold for all other cycles $C$.
Thus we have have $2^{\nu-1}$ choices for $(\boldsymbol{\mu}^e,\boldsymbol{\mu}^o)$. 
\si

For each of those choices, (4.28) becomes
$$
\Cal K_{(\boldsymbol{\mu}^e,\boldsymbol{\mu}^o)}(z)
=\prod_{i\in [\nu]}z^{3\mu_i}
\,\cdot\prod_{1\le i<j\le \nu\atop \mu_{ij}^e=\mu_{ij}}
\!\!\!\!z^{2(\mu_{ij}-1)}\!\!\left[\mu_{ij}-(\mu_{ij}-1)z^2\right]
\,\cdot\prod_{1\le i<j\le \nu\atop \mu_{ij}^e=0}\!\!\!
z^{\mu_{ij}}.
\tag 4.29
$$
For each of these $2^{\nu-1}$ polynomials,
$$
\Cal K_{(\boldsymbol{\mu}^e,\boldsymbol{\mu}^o)}(1) =1,\tag 4.30
$$
and
$$
\text{deg}\,\Cal K_{(\boldsymbol{\mu}^e,\boldsymbol{\mu}^o)}(z)\le 3\sum_i\mu_i+2\sum
_{1\le i<j\le \nu}\mu_{ij}.
$$ 
Using the constraints
$$
\sum_i\mu_i+\sum_{1\le i<j\le \nu}\mu_{ij}=\mu,\quad 
\sum_i\mu_i+2\sum_{1\le i<j\le \nu}\mu_{ij}\ge 3\nu,
$$
we easily obtain then that
$$
\text{deg}\,\Cal K_{(\boldsymbol{\mu}^e,\boldsymbol{\mu}^o)}(z)\le 4\mu-3\nu.
$$

Now, the generating function $\Cal H_M^e(z)$ of the core-type graphs $G$ without
odd cycles that cancel to $M$ is the sum of 
$\Cal H_{(\boldsymbol{\mu}^e,\boldsymbol{\mu}^o)}(z)$ over all $2^{\nu-1}$ sets of
feasible pairs $(\boldsymbol{\mu}^e,\boldsymbol{\mu}^o)$. Using (4.27), (4.29) and (4.30)
we arrive at the following formula.
\proclaim{Lemma 4.4} For each kernel $M$, with $\ell:=\mu-\nu>0$,
$$
\Cal H_M^e(z)=\frac{\kappa}{\nu!}\frac{2^{\nu-1}}{(1-z^2)^{\mu}}P_M(z),\tag 4.31
$$
where $P_M(z)$ is a polynomial of degree $4\mu-3\nu=\mu+3\ell$ at most, and $P_M(1)=1$.
\endproclaim
This Lemma directly implies 
\proclaim{Corollary 4.5}
$$
C_{\ell}^e(x)=\sum_{M:e(M)-v(M)=\ell}H_M^e(x),\tag 4.32
$$
where
$$
H_M^e(x)=\frac{\kappa}{\nu!}\frac{2^{\nu-1}}{(1-T^2(x))^{\mu}}P_M(T(x)),\tag 4.33
$$
\endproclaim

Using $\ell<\mu(M)\le 3\ell$, $\nu(M)=\mu(M)-\ell$, we deduce from
(4.32)-(4.33) that
$$
C_{\ell}^e(x)=(1-T^2(x))^{-3\ell}\left[\sum_{d=0}^{8\ell-1}c_{\ell,d}^e(1-T(x))^d
\right],\tag 4.34
$$
where
$$
c_{\ell,0}^e=\frac{2^{2\ell-1}}{(2\ell)!}\sum_{\boldsymbol\mu\text{ meets }(4.19)}
\kappa(\boldsymbol\mu).\tag 4.35
$$
So, by the second line in (4.22), $c_{\ell,0}^e=2^{2\ell-1}c_{\ell,0}$. The
proof of Proposition 4.3 is complete.\qed
\bi

Comparing (4.14) and (4.34)-(4.35), and using $T(e^{-1})=1$, we obtain: for $\ell>0$,
$$
C_{\ell}^e(x)=\frac{1}{2^{\ell+1}}C_{\ell}(x)+O\big(|1-T(x)|^{-3\ell+1}\bigr),\quad
(|x|<e^{-1},\,x\to e^{-1}).\tag 4.36
$$
And we recall, (4.3), that
$$
C_0^e(x)=\frac{1}{2}C_0(x)+\ln(2^{-1/4}e^{1/8})+O(|T(x)-1|).\tag 4.37
$$
Now, by (4.1), for a fixed $L>0$,
$$
\pr(S_n>0,\,\Cal E_n\le L)=N(n,p)\oint_{\Gamma}x^{-n-1}\exp\left[\sum_{\ell=-1}^L
\left(\frac{p}{q}\right)^{\ell}C_{\ell}^e(x)\right]\,dx,
$$
where $\Gamma$ is within the disc $|x|<e^{-1}$. As in Section 2.3, we switch
to $y$ by $x=ye^{-y}$, and choose in the $y$-plane the circular contour $\Gamma^\prime$ 
$y=e^{-an^{-1/3}+i\theta}$, $a>0$ being fixed this time. Observe that, 
for each $1\le\ell\le L$,
$$
\left(\frac{p}{q}\right)^{\ell}|1-T(ye^{-y})|^{-3\ell+1}
\le_b n^{-\ell}|1-y|^{-3\ell+1}\le_b n^{-1/3},
$$
and, likewise in (4.3) the remainder term $O(|T(ye^{-y})-1|)$ is 
$O(n^{-1/3})$. And
of course $C_{-1}^e(ye^{-y})=C_{-1}(ye^{-y})$. On the basis of (4.36)-(4.37), it
can be shown then that
$$
\multline
\pr(S_n>0,\,\Cal E_n\le L)\\
\sim 2^{-1/4}e^{1/8}N(n,p)\oint_{\Gamma^\prime}(ye^{-y})^{-n-1}
\exp\left[\frac{1}{2}\sum_{\ell=-1}^L
\left(\frac{p}{2q}\right)^{\ell}C_{\ell}(ye^{-y})\right]\,d(ye^{-y}),
\endmultline
$$
where now $\Gamma^\prime$ can be replaced by a circular contour of an
arbitrarily small radius. Going back to the $x$-plane, we recognize (see (2.1.11)) 
the value of
the resulting integral as 
$$
2^{-1/4}e^{1/8}\left.\pr(S_n>0,\,\Cal E_n\le L)\right|_{\hat p=1/2},
$$
for $\la=O(1)$ needless to say. By (2.1.13) the latter probability is at least
$$
\left.\pr(S_n>0)\right|_{\hat p=1/2}\,\cdot\pr(\Cal E_n\le L).
$$ 
Letting $n\to\infty$, and using the part {\bf (i)\/} for
$\left.\pr(S_n>0)\right|_{\hat p=1/2}$, we get
$$
\liminf\frac{\pr(S_n>0,\,\Cal E_n\le L)}{2^{-1/4}e^{1/8}c(\la)n^{-1/12}}\ge 
\liminf \pr(\Cal E_{n}\le L).
$$
Since $\pr(S_n>0)\ge \pr(S_n>0,\,\Cal E_n\le L)$, and $\Cal E_n=O_P(1)$, letting
$L\uparrow\infty$ enables us to conclude that
$$
\pr(S_n>0)\gtrsim 2^{-1/4}e^{1/8}c(\la)n^{-1/12}.
$$
Together with (4.14) this proves that
$$
\pr(S_n>0)\sim 2^{-1/4}e^{1/8}c(\la)n^{-1/12}.
$$
The proof of Theorem 1.1 {\bf (ii)\/} is now complete.\qed
\bi

{\bf Acknowledgement.\/} We are grateful to the participants of a 
graduate student workshop at Ohio State University for helpful commentaries
during numerous discussions of various phases of this study. We thank 
Mike Molloy for posing a $2$-colorability problem of the critical graph 
$G(n,m)$, and Greg Sorkin for an encouragement when we needed it most.
\'Akos Seres and Saleh Tanveer, members of the second author's  Ph.D Committee,
provided valuable comments on the project.
\bi

\bi
\si
{\bf Appendix.\/}
\si
{\bf Proof of (2.3.22).\/} {\bf (i)\/} For the second bound, we use (2.3.19) and, setting
$s=a+it$,
$$
\left|s^{1-y}\exp\left(\frac{\mu s^2}{2}+\frac{s^3}{3}\right)\right|\le
a^{1-y}\exp\left(\frac{\mu a^2}{2}+\frac{a^3}{3}\right)
\exp\bigl[-t^2(a+\mu/2)\bigr].
$$
So
$$
A(y,\mu)\le (2\pi)^{-1}\sqrt{\frac{\pi}{a+\mu/2}}\exp\left(-
\frac{\mu^3}{6}+\frac{\mu a^2}{2}+\frac{a^3}{3}\right).
$$
\qed
\si

{\bf (ii)\/} For the first bound, we use (2.3.20), i. e.
$$
e^{\mu^3/6}\,3^{(y+1)/3}A(y,\mu)=\sum_{k=0}^{\infty}
\frac{(3^{2/3}\mu/2)^k}{\Gamma(k+1)\Gamma((y+1-2k)/3)},\tag A.1
$$
and the inequalities
$$
\frac{(1/2)^{a+b+2}}{a+b+2}\le\frac{\Gamma(a+1)\Gamma(b+1)}{\Gamma(a+b+2)}
\le \frac{a^ab^b}{(a+b)^{a+b}}\le 1,\quad (a\ge 0,\,b\ge 0),
$$
which follow from a classic formula
$$
\int_0^1x^a(1-x)^b\,dx =\frac{\Gamma(a+1)\Gamma(b+1)}{\Gamma(a+b+2)},
$$
and 
$$
\align
\max_{x\in [0,1]}x^a(1-x)^b=&\,\frac{a^ab^b}{(a+b)^{a+b}},\\
\int_0^1x^a(1-x)^b\,dx \ge&\, (1/2)^b\int_0^{1/2}x^a\,dx+(1/2)^a\int_{1/2}^1
(1-x)^b\,dx\\
=&\,(1/2)^{a+b+1}\left(\frac{1}{a+1}+\frac{1}{b+1}\right).
\endalign
$$

Break the sum in (A.1) into $\Sigma_1$, $\Sigma_2$, and $\Sigma_3$, for
$\{k\ge 2:\,(y+1-2k)/3\ge 1\}$, $\{k\ge 1:\,(y+1-2k)/3\le 0\}$, and $\{k=0,1:\,
\text{ or }(y+1-2k)/3\ge 1\}$, respectively.
(Recall that $\Gamma(0)=\infty$.) For $\Sigma_1$,
$$
\align
\frac{1}{\Gamma(k+1)\Gamma((y+1-2k)/3))}=&\,\frac{\Gamma(2k/3)}{\Gamma(k+1)}
\cdot\frac{1}{\Gamma((y+1-2k)/3)\Gamma(2k/3)}\\
\le&\,\frac{2^{(y+1)/3}}{\Gamma((y+1)/3)}\cdot\frac{\Gamma(2k/3)\Gamma(k/3+1)}
{\Gamma(k+1)}\cdot\frac{1}{\Gamma(k/3+1)}\\
\le&\,\frac{2^{(y+1)/3}}{\Gamma((y+1)/3)}\cdot 
2\,\frac{(2k/3)^{2k/3-1}(k/3)^{k/3}}{k^{k-1}} \cdot\frac{1}{\Gamma(k/3+1)}\\
\le&\,6\,\frac{2^{(y+1)/3}}{\Gamma((y+1)/3)}\cdot\frac{(2^{2/3}/3)^{k}}
{\Gamma(k/3+1)}.
\endalign
$$
Therefore
$$
|\Sigma_1|\le 6\,\frac{2^{(y+1)/3}}{\Gamma((y+1)/3}\sum_{k\ge 0}\frac
{(|\mu|^3/6)^{k/3}}{\Gamma(k/3+1)}\le_b(|\mu|^3\lor 1)
\frac{2^{(y+1)/3}}{\Gamma((y+1)/3)}\cdot e^{|\mu|^3/6}.\tag A.2
$$
For $\Sigma_2$, we use 
$$
\Gamma(z)\Gamma(1-z)=\frac{\pi}{\sin(\pi z)}\Longrightarrow\frac{1}{|\Gamma((y+1-2k)/3)|}
\le\Gamma(1+(2k-y-1)/3),
$$
and
$$
\align
&\frac{\Gamma(1+(2k-y-1)/3)}{\Gamma(k+1)}
\le\,\frac{1}{\Gamma(1+(y+1)/3)}\cdot\frac{\Gamma(2k/3+2)}{\Gamma(k+1)}\\
=&\,\frac{1}{\Gamma(1+(y+1)/3)}
\cdot\frac{\Gamma(2k/3+2)\Gamma(k/3+1)}{\Gamma(k+1)}\cdot\frac{1}{\Gamma(k/3+1)}\\
\le&\,\frac{1}{\Gamma(1+(y+1)/3)}\cdot\frac{\Gamma(k+3)}{\Gamma(k+1)}\cdot\frac{(k/3)^{k/3}(1+2k/3)^{1+2k/3}}{(k+1)^{k+1}}\\
\le&\,\frac{1}{\Gamma(1+(y+1)/3)}\cdot (k+2)^2\frac{(2^{2/3}/3)^k}{\Gamma(k/3+1)}.
\endalign
$$
Therefore
$$
|\Sigma_2|\le_b(|\mu|^9\lor 1)\frac{1}{\Gamma((y+1)/3)}\cdot e^{|\mu|^3/6}.\tag A.3
$$
And it is not difficult to show that
$$
|\Sigma_3|\le_b |\Sigma_1|+|\Sigma_2|.\tag A.4
$$
The relations (A.1)-(A.4) imply that
$$
A(y,\mu)\le_be^{|\mu|^3/2}\,\frac{(2/3)^{\frac{y+1}{3}}}{\Gamma\left(\frac{y+1}{3}\right)}.
$$
\qed
\bi
{\bf Proof of (2.3.23)-(2.3.24).\/} Again we use (2.3.19). Let us choose $a=\xi$, where $\xi
=\xi(y,\mu)$ is a maximum point of
$$
\Psi(a;y,\mu):=-y\ln a+\frac{\mu a^2}{2}+\frac{a^3}{3},\quad a\in (0,\infty),
$$
i. e. a positive root of
$$
\Psi^{(1)}_a(a;y,\mu)=\mu a +a^2-\frac{y}{a}=0.\tag A.5
$$
A root  exists and is unique, since $\Psi_a(0+;y,\mu)=-\infty$, $\Psi(\infty;y,\mu)=\infty$ and
$$
\Psi^{(2)}_a(a;y,\mu)=\mu+2a+\frac{y}{a^2}>0,\quad (a\ge 0).
$$
Observe that $\mu\xi^2/y$ is bounded away from zero. If not, then,  by (A.5),
$$
\frac{\mu^3\xi^6}{y^3}\to 0,\,\,\frac{y^2}{\xi^6}\to 1,
$$
which implies that $\mu^3/y\to 0$, contradicting $y=O(\la^3)=O(\mu^3)$.
\si
Break the integral in (2.3.19) into $I_1$ over $|t|\le \mu^{-1/2}y^{1/7}$ 
and $I_2$ over $|t|\ge
\mu^{-1/2}y^{1/7}$. Arguing as the part {\bf (i)\/} of the previous proof, 
we bound
$$
\aligned
|I_2|\le_b&\,\xi\exp[\Psi(\xi;y,\mu)]\!\!\!\int\limits_{|t|\ge 
\mu^{-1/2}y^{1/7}}\!\!\!\!\!\!\!
\exp\left[-t^2(\xi+\mu/2)\right]\,dt\\
\le_b&\,\frac{\xi\exp[\Psi(\xi;y,\mu)]}{(\xi+\mu/2)^{1/2}}\cdot 
e^{-y^{2/7}/2}.
\endaligned\tag A.5
$$
Turn to $I_1$. Since 
$$
\Psi_s^{(3)}(s;y,\mu)=-\frac{2y}{s^3}=O(y\xi^{-3}),
$$
we have
$$
\Psi(s;y,\mu)=\Psi(\xi;y,\mu)-\frac{t^2}{2}(\mu+2\xi+y\xi^{-2})
+O\bigl(y\xi^{-3}\mu^{-3/2}y^{3/7}\bigr),
$$
and
$$
y\xi^{-3}\mu^{-3/2}y^{3/7}=\frac{y^{3/7-1/2}}{\left(\frac{\mu\xi^2}{y}
\right)^{3/2}}\le_b y^{-1/14}.
$$
Consequently
$$
\aligned
I_1\sim&\,\xi\exp[\Psi(\xi;y,\mu)]\int\limits_{|t|\le 
\mu^{-1/2}y^{1/7}}\exp\left(-\frac{t^2}{2}(\mu+2\xi+y\xi^{-2})\right)\,
dt\\
\sim&\,\xi\exp[\Psi(\xi;y,\mu)]\left(\frac{2\pi}{\mu+2\xi+y\xi^{-2}}
\right)^{1/2}.
\endaligned\tag A.6
$$
Since $y\xi^{-2}=O(\mu)$, (A.5)-(A.6) imply that $I_1\gg I_2$, hence
$$
A(y,\mu)\sim e^{-\mu^3/6}(2\pi)^{-1}I_1,
$$
which proves (2.3.23). \qed
\bi

If we drop the condition $y\to\infty$, then the integral in (2.3.19)
is of order
$$
\xi\exp[\Psi(\xi;y,\mu)]\int\limits_{-\infty}^{\infty} 
\exp\left[-t^2(\xi+\mu/2)\right]\,dt
=\xi\exp[\Psi(\xi;y,\mu)]\left(\frac{2\pi}{\xi+\mu/2}\right)^{1/2},
$$
which proves (2.3.24). \qed
\bi
{\bf Proof of (1.5).\/} The system (1.1) is solvable iff for every cycle $C$ of $G$,
$$
\sum_{e\in E(C)}b_e=O(\text{ mod }2).\tag A.7
$$
If $b_e\in\{0,1\}$ are independent random variables with $\pr(b_e=1)=\hat p$,
the condition (A.7) is met with probability $(1+(1-2\hat p)^{|C|})/2$, Kolchin [13].
\si

Consider $G=G(n,p=\ga/n)$, $\ga<1$. Let $X_{ns}$ denote the number
of cycles of length $s$ which are ``bad'', i. e. do not meet the condition (A.7).
We need to find the limiting 
distribution of $X_n=\sum_{s\ge 3}X_{ns}$ the total number of  ``bad'' cycles.
To this end, observe that, with probability approaching $1$,
the cycles $G(n,p)$ may have are those in the unicyclic components. Let us call
them u-cycles. The expected
number of all cycles of length $k\ge 3$ is
$$
\binom{n}{k}\frac{(k-1)!}{2}p^k\le \frac{\ga^k}{2k}.
$$
So
$$
\lim_{A\to\infty}\lim_{n\to\infty}\pr(G(n,p)\text{ has a cycle of length }\ge A)=0.
$$
 Let $Y_{ns}$ be the total number of all u-cycles
of length $s$. In [19] it was proven that, for $\ga$ fixed, $\{Y_{ns}\}_{
s\le A}$ converges in distribution to $\{\text{Poisson}(\sigma_s)\}_{s\le A}$,
where the Poissons are independent and
$$
\sigma_s=\frac{T^s(\ga e^{-\ga})}{2s},\quad s\ge 3.
$$
As $\ga<1$, we have $T(\ga e^{-\ga})=\ga$, because $T(x)=xe^{T(x)}$, for
$x<e^{-1}$. Now a u-cycle of length $s$ is bad with probability
$$
\pi_s=\frac{1-(1-2\hat p)^s}{2}.
$$
Consequently $\{X_{ns}\}_{s\le A}$ converges to 
$\{\text{Poisson}(\pi_s\sigma_s)\}_{s\le A}$, whence $X_n$ converges to
$\text{Poisson }\left(\sum_{s\le A}\pi_s\sigma_s\right)$.
Therefore, 
$$
\lim_{n\to\infty}\pr\{\text{there are nod bad u-cycles of length }A
\text{ at most}\}
=e^{-\sum_{s\le A}\pi_s\sigma_s}.
$$
It remains to notice that
$$
\align
\sum_{s\ge 3}\pi_s\sigma_s=&\,\sum_{s\ge 3}\frac{1-(1-2\hat p)^s}{2}\frac{\ga^s}{2s}\\
=&\,\frac{1}{4}\ln\frac{1-\ga(1-2\hat p)}{1-\ga}-\frac{\ga}{2}\hat p-
\frac{\ga^2}{2}\hat p(1-\hat p).
\endalign
$$
That the same formula holds for $G(n,m=\ga n/2)$ follows then in a standard way.
\qed

\enddocument